# A Bounded Arithmetic $AID$ for Frege System

Toshiyasu Arai
Faculty of Integrated Arts and Sciences
Hiroshima University
Higashi-Hiroshima, 739-8521 Japan
e-mail: arai@mis.hiroshima-u.ac.jp
and
The Fields Institute
Toronto, Ontario, M5T3J1 Canada

Apr. 3, 1998

**Abstract**

In this paper we introduce a system $AID$ (Alogtime Inductive Definitions) of bounded arithmetic. The main feature of $AID$ is to allow a form of inductive definitions, which was extracted from Buss' propositional consistency proof of Frege systems $\mathcal{F}$ in [7]. We show that $AID$ proves the soundness of $\mathcal{F}$, and conversely any $\Sigma_0^b$-theorem in $AID$ yields boolean sentences of which $\mathcal{F}$ has polysize proofs. Further we define $\Sigma_1^b$-faithful interpretations between $AID + \Sigma_0^b - CA$ and a quantified theory $QALV$ of an equational system $ALV$ in P. Clote [10]. Hence $ALV$ also proves the soundness of $\mathcal{F}$.

There are two sources by Cook [14] and Buss [7] to motivate this paper.

In the pioneering paper [14] Cook shows that his equational theory $PV$ corresponds to the *extended Frege system* $e\mathcal{F}$. This means that $PV$ proves the soundness of $e\mathcal{F}$, and any provable equation in $PV$ can be transformed into boolean tautologies so that $e\mathcal{F}$ has polysize proofs of these tautologies. Thus extended Frege system has also polysize proofs of partial consistencies of himself. Later Buss [4] shows that the same thing holds for the $\Sigma_1^b$-theorems of the bounded arithmetic $S_2^1$ in place of $PV$ and $e\mathcal{F}$. Note that these theories $PV$ and $S_2^1$ characterize the complexity class $P$.

It has remained open what bounded arithmetic $T_\mathcal{F}$ corresponds to the *Frege system* $\mathcal{F}$. That is, a bounded arithmetic $T_\mathcal{F}$ such that $T_\mathcal{F}$ proves the soundness of $\mathcal{F}$, and for any formula in a restricted form if $T_\mathcal{F}$ proves the formula, then it can be transformed into boolean sentences so that $\mathcal{F}$ has polysize proofs of these sentences.

In [7] Buss shows an intensionally correct truth definition for boolean sentences can be written as polysize boolean formulae. Thus Frege system has also polysize proofs of partial consistencies of himself. The truth definition utilizes countings, vector summations and a form of inductive definitions. Prior to this Buss [5] shows that Frege system $\mathcal{F}$ has polysize proofs of the propositional pigeon hole principles by showing that $\mathcal{F}$ has an intensionally correct definitions of counting. The definition utilizes carry-save-additions. Therefore some nontrivial parts of mathematics are carried in the propositional proof system $\mathcal{F}$. Such a fact has been already found by Cook [14] for extended Frege system $e\mathcal{F}$. A natural question to be asked is: Is this all for $\mathcal{F}$? Namely all parts of mathematics included in $\mathcal{F}$ are derived from counting?

The former ones, viz. countings and vector summations suffices to develop some metamathematics-arithmetizations in $\mathcal{F}$. It seems to us that the latter ingredient, viz. inductive definition is essential. In fact the former are derived from the latter. This is shown in Section 2 and is expected: The latter inductive definition corresponds to the evaluation for the computation tree of a predicate in the complexity class $ALOGTIME$. Prior to [7] Buss [6] shows that BSVP (Boolean Sentence Value Problem) is in $ALOGTIME$ and hence is $ALOGTIME$-complete in a weak reducibility. In view of this result in [6] the task in [7] is to show that the algorithm for $BSVP \in ALOGTIME$ is intensionally correct for $\mathcal{F}$.

In this paper $ALOGTIME$ is used as a synonym of (uniform) $NC^1$. By definition a function $f$ of polynomial growth rate is $ALOGTIME$-computable, denoted $f \in \mathcal{F}ALOGTIME$, iff its $bitgraph$ $\{i < |f(\bar{x})| : Bit(i, f(\bar{x})) = 1\}$ is in $ALOGTIME$.

When one reads these proofs in $\mathcal{F}$ in [5], [7], it is natural to consider that these are images of mathematical-logical-arithmetical proofs in a system of bounded arithemetic. There may be various ways to formulate preimages, i.e., a system of bounded arithmetic in which proofs proofs in [5], [7] are carried out. In this paper we introduce a system $AID$ (Alogtime Inductive Definitions) of bounded arithmetic. The main feature of $AID$ is to allow a form of inductive definitions. We show the following results.

1. (*cf.* Theorem 2.1.) Bounded vector summation and hence Counting are $\Sigma_0^b$-bitdefinable in $AID$ (carry-save-addition reduces to inductive definitions in $AID$): If the bitgraph of an $f(i, \bar{y})$ is $\Sigma_0^b$-definable in $AID$, then so is the function
$$g(x, \bar{y}) = \sum_{i < |x|} f(i, \bar{y})$$

2. (*cf.* section 4.) In $AID$ a truth definition for boolean sentences is definable by a $\Sigma_0^b$-formula $TRUE(x)$. The definition is nothing but the arithmetical equivalent to Buss' definition in [7] and hence is intensionally correct. Therefore

3. (*cf.* Theorem 4.1.) $AID \vdash RFN(\mathcal{F})$, where $RFN(\mathcal{F})$ denotes the reflection schema for $\mathcal{F}$.

4. (*cf.* section 6.) For each $\Sigma_0^b$-formula $B(x)$ there exists a $\Sigma_0^b$-bitdefinable function $\varphi_B : x \mapsto \langle B(x) \rangle$ such that $\langle B(x) \rangle$ is a boolean sentence and $AID \vdash B(x) \leftrightarrow TRUE(\langle B(x) \rangle)$.

5. (*cf.* Theorem th:AIDF.) If $AID \vdash B(x)$ for a $\Sigma_0^b$-formula $B(x)$, then $AID + \Sigma_0^b - CA \vdash \mathcal{F} \vdash^{p|x|} \langle B(x) \rangle$ for a plynomial $p \, |x|$. Thus

$$AID \vdash B(x) \Leftrightarrow AID + \Sigma_0^b - CA \vdash \mathcal{F} \vdash^{p|x|} \langle B(x) \rangle \text{ for a plynomial } p \, |x|$$

6. (*cf.* Theorem 8.2.) A predicate is in $ALOGTIME$ iff it is $\Sigma_0^b$-definable in $AID$ iff it is $\Delta_1^b$-definable in $AID$.

The aim of $AID$ is to capture, calibrate and draw the line to mathematical-arithmetical power of Frege system $\mathcal{F}$. By establishing what can be done in $\mathcal{F}$ mathematically, we hope to find a right candidate of hard tautologies for $\mathcal{F}$ *cf.* [3] and [16], to specify what kind of nonstandard models to be considered in order to prove superpolynomial lowerbounds for hard tautologies for $\mathcal{F}$, *cf.* $I\Delta_0(f)$ vs. constant depth Frege and nonstandard model of $PA$ in Ajtai [1], and to find a bounded arithmetic for constant depth threshold gates $TC^0$ in order to compare $AID$.

Let us explain contents of sections.

In section 1 the bounded arithmetic $AID$ is defined. First we introduce a base language $L_{BA}$ of weak bounded arithmetics and then the language of $AID$ is obtained by adding predicate constants for inductively defined predicates. Also we recall some axiom schemata.

In section 2 we show first that inductive definitions along quadtree (4-branching tree), iterated inductive definitions and simultaneous inductive definitions reduce to a single inductive definition available in $AID$. Using these we show second that vector summations are $\Sigma_0^b$-bitdefinable in $AID$.

In section 3 we show that each predicate in $ALOGTIME$ is $\Sigma_0^b$-definable in $AID$.

In section 4 we examine Buss' propositional consistency proof of Frege systems in [7] and verify that it is formalizable in $AID$.

In section 5 *stratifications* of formulae are defined. These are in essence to interpret first order formulae into second order formulae. In later sections we need these.

In section 6 we show that any $\Sigma_0^b$-theorem in $AID$ yields true boolean sentences of which $\mathcal{F}$ has polysize proofs.

In section 7 we introduce some systems of bounded arithmetic in the language $L_{BA}$, i.e., without inductively defined predicate which are equivalent to $AID$.

In section 8 we show that $\Sigma_1^b$-consequences in, e.g., $AID + \Sigma_0^b - CA$ can be realized by a $\Sigma_0^b$-set $\{i < p \mid \bar{x} \mid : A(\bar{x}, i)\} \bar{x})$.

In [9] P. Clote defines a function algebra $N_0$ and show that $N_0 = \mathcal{F}ALOGTIME$, the class of $ALOGTIME$-computable functions. Then he introduces an equational system $ALV$ based on $N_0$ in [10]. In section 9 we show that $AID + \Sigma_0^b - CA$ is equivalent to $ALV$ in the sense that there exist $\Sigma_1^b$-faithful interpretations between $AID + \Sigma_0^b - CA$ and a quantified theory $QALV$. Hence we see that Clote's $ALV$ proves the soundness of Frege system, Corollary 9.1.

There are now several theories besides $ALV$ which are designed for Frege systems in Clote [11] and for $ALOGTIME$ in Clote and Takeuti [12], [13], and F. Pitt [18]. It is open to us whether these are equivalent to $AID$.

The results in section 1-8 of this paper were contained in a handwritten note 'Frege System, $ALOGTIME$ and Bounded Arithmetic' written in Dec. 1991. Section 9 is augmented here.

**Acknowledgements**. I would like to thank S. A. Cook for inviting me to Fields Institute, Toronto, hospitality during my stay and encouraging me to complete the work on $AID$. Without his interests and encouragements it would be impossible to finish this paper.

# 1 The bounded arithmetic $AID$

In this section we introduce the bounded arithmetic $AID$.

## 1.1 A base language

Function constants in the *language $L_{BA}$* of weak bounded arithmetic are 0 (zero), 1 (one), + (addition), $\lfloor \frac{x}{2} \rfloor$ (half), $|x|$ (length), $x \# y$ (smash), $x \cdot 2^{|y|}$ (padding), $x \dot{-} y$ (modified subtraction) and $x[i, j]$ (part). Let $x_k$ denotes the $k$-th digit of $x$ in binary representation. Then the part function $x[i, j]$ is defined to be the string $x_{j-1} \cdots x_i$ from $i$-th digit $x_i$ to $j - 1$-th digit $x_{j-1}$:

$$x[i, j] = \sum_{i \leq k < j} x_k \cdot 2^{k-i} \text{ for } x = \sum_k x_k \cdot 2^k$$

Clearly
$$\mid x[i, j] \mid = \min\{j, \mid x \mid\} \dot{-} i$$

Thus $L_{BA}$ is obtained from the language of $S_2^i$ in [4] by deleting the multiplication $x \cdot y$ and adding three functions $x \cdot 2^{|y|}, x \dot{-} y, x[i, j]$. Further the successor function $Sx$ is replaced by $x + 1$. From these familiar functions are defined as follows:

1. $\mid x \mid_0 = x$ and $\mid x \mid_{n+1} = \mid\mid x \mid_n \mid$. Also $\|x\| = \mid x \mid_2$.

2. $Bit(i, x) = x[i, i+1)$ ($i$-th digit)

3. $MSP(x, i) x[i, |x|)$ (Most Significant Part)

4. $|x| \cdot |y| = |x \# y| \dot{-} 1$. *Multiplication for small numbers* is definable as follows:
$$multi(i, j, x, y) = \begin{cases} i \cdot j = |x[0, i)| \cdot |y[0, j)| & \text{if } i \leq |x| \ \& \ j \leq |y| \\ 0 & \text{otherwise} \end{cases}$$

5. *Successor functions $xi$ in binary notation*: Put
$$xi =_{df} s_i x = 2 \cdot x + i = x + x + i \text{ for } i < 2.$$

We encode a word $i_{k-1} \cdots i_0 \in \{0,1\}^*$ by *attaching the leading marker* 1:
$$\lceil i_{k-1} \cdots i_0 \rceil = 1 i_{k-1} \cdots i_0 \text{ (in binary notation)} \tag{1}$$

For example $\lceil \varepsilon \rceil = 1$ for the empty word $\varepsilon$. Using this encoding *concatenation* on words is defined as follows:
$$x * y = x \cdot 2^{|y| \dot{-} 1} + y[0, |y| \dot{-} 1)$$

If $x = \lceil i_{k-1} \cdots i_0 \rceil$, $y = \lceil j_{l-1} \cdots j_0 \rceil$, then $x * y = \lceil i_{k-1} \cdots i_0 j_{l-1} \cdots j_0 \rceil$. Cleraly $|x * y| = |x| + |y| \dot{-} 1$ for $y > 0$.

$BASIC$ denotes the set of basic axioms for constants in $L_{BA}$. These are obtained from basic axioms in p.31, [4] by deleting axioms mentioning multiplication and adding the following axioms:

1. $j \leq i \to x[i, j) = 0$ ; $i < j \to x[i+1, j+1) = \lfloor \frac{2}{x} \rfloor [i, j)$

2. $(2x)[0, j+1) = 2 \cdot (x[0, j)) \ \& \ (2x+1)[0, j+1) = 2 \cdot (x[0, j)) + 1$

3. $x \cdot 2^{|0|} = x$ ; $y \neq 0 \to x \cdot 2^{|y|} = x \cdot 2^{|\lfloor \frac{2}{y} \rfloor|} + x \cdot 2^{|\lfloor \frac{2}{y} \rfloor|}$

4. $x \leq y \to x \dot{-} y = 0$ ; $x \geq y \to (x+1) \dot{-} y = (x \dot{-} y) + 1$

When the language is expanded to include a set $\mathcal{X}$ of predicate constants, then $BASIC$ is assumed to include the equality axiom for constants in $\mathcal{X}$. The expanded language is denoted $L_{BA}(\mathcal{X})$.

Quantifiers of the form $Qx \leq t$ ($Q \in \{\forall, \exists\}$) are said to be *bounded quantifiers* while quantifiers of the form $Qx \leq |t|$ ($Q \in \{\forall, \exists\}$) are *sharply bounded quantifiers*.

Classes of $\Sigma_i^b$ formulae and $\Pi_i^b$ formulae are defined as in [4]. Also a formula $A$ is in $s\Sigma_i^b$ (*strict* $\Sigma_i^b$) iff $A$ is in a prenex normal form such that its leading quantifier is a bounded existential quantifier followed by a string of alternating bounded quantifiers with a sharply bounded matrix $B \in \Sigma_0^b$:
$$A \equiv \exists x_1 \leq t_1 \forall x_2 \leq t_2 \cdots Q x_i \leq t_i B$$

$s\Pi_i^b$ is defined dually. Relativzed classes $\Sigma_i^b(\mathcal{X}), \Pi_i^b(\mathcal{X}), s\Sigma_i^b(\mathcal{X})$ and $s\Pi_i^b(\mathcal{X})$ are defined analogously. $\Sigma_i^b(\mathcal{X})$, etc. is denoted $\Sigma_i^b(L)$ for $L = L_{BA}(\mathcal{X})$.

$\Sigma_0^b - LIND$ denotes the following axiom schema:

$$A(0) \wedge \forall y < |x| \, (A(y) \to A(y+1)) \to A(|x|)$$

for $A \in \Sigma_0^b$.

A base fragment $\Sigma_0^b - LIND$ of bounded arithmetic: its language is $L_{BA}$ and its axioms are $BASIC$, the axiom schema $\Sigma_0^b - LIND$ and the *Bit Extensionality Axiom*:

$$|x| = |y| \, \& \, \forall i < |x| \, (Bit(i, x) = Bit(i, y)) \to x = y$$

A term of the form $\ell\bar{x} = \sum_{i=1}^{n} c_i x_i + d$ for some constants, i.e, numerals $c_i \, (1 \leq i \leq n), d$ is said to be a *linear form* (in a sequence $\bar{x} = x_1, \ldots, x_n$ of variables). Also $\ell\|\bar{x}\|$ denotes $\sum_{i=1}^{n} c_i \|x_i\| + d$.

Let $t$ be either a polynomial $p \, |\bar{x}| = p(|x_1|, \ldots, |x_n|)$ or a linear form $\ell\|\bar{x}\|$. Assume that a variable $y$ does not occur in $t$. Then $\forall \, |y| \leq t \, B(y)$ denotes the formula $\forall y \leq 2^t (|y| \leq t \to B(y))$. If $t$ is a polynomial $p \, |\bar{x}|$, then the quantifier $\forall \, |y| \leq t$ is a bounded quantifier since the smash function $\#$ is in $L_{BA}$. If $t$ is a linear form, then it is a sharply bounded quantifier. The existential one $\exists \, |y| \leq t$ is defined dually.

## 1.2 Some axiom schemata

Put

$$i \in x \Leftrightarrow_{df} Bit(i, x) = 1$$

Using this we define analogues of some axiom schemata in second order arithmetic. Let $L_{BA}(\mathcal{X})$ be an expanded language for a set $\mathcal{X}$ of predicate constants and $\Phi$ a set of formulae in $L_{BA}(\mathcal{X})$.

1. $\Phi - CA$ denotes the axiom schema:

   $$\exists \, |y| \leq p \, |\bar{x}| \, \forall i < p \, |\bar{x}| \, (i \in y \leftrightarrow B(i; \bar{x}))$$

   for each polynomial $p \, |\bar{x}|$ and each formula $B(i; \bar{x}) \in \Phi$.

2. $\Phi - LCA$ denotes the axiom schema:

   $$\exists \, |y| \leq \ell\|\bar{x}\| \forall i < \ell\bar{x} \, (i \in y \leftrightarrow B(i; \bar{x}))$$

   for each linear form $\ell\|\bar{x}\|$ and each formula $B(i; \bar{x}) \in \Phi$.

3. $\Delta_1^b(L_{BA}(\mathcal{X})) - CA$ denotes the axiom schema:

   $$\forall i (B(i; \bar{x}) \leftrightarrow \neg C(i; \bar{x})) \to \exists \, |y| \leq p \, |\bar{x}| \, \forall i < p \, |\bar{x}| \, (i \in y \leftrightarrow B(i; \bar{x}))$$

   for each polynomial $p \, |\bar{x}|$ and each $B(i; \bar{x}), C(i; \bar{x}) \in \Sigma_1^b(L_{BA}(\mathcal{X}))$.

4. $\Phi - AC$ (or $\Phi$-replacement) denotes the axiom schema:

$$\forall i < p\,|\bar{x}|\,\exists\,|y| \leq q\,|\bar{x}|\,B(i,y;\bar{x}) \to$$
$$\exists\,|z| \leq p\,|\bar{x}|\cdot q\,|\bar{x}|\,\forall i < p\,|\bar{x}|\,B(i,z_i;\bar{x})$$

for polynomials $p\,|\bar{x}|, q\,|\bar{x}|$ and each $B \in \Phi$, where $z_i = z[q\,|\bar{x}|\cdot i, q\,|\bar{x}|\cdot(i+1))$.

Let $\Delta_1^b(L_{BA}(\mathcal{X})) - LIND$ denotes the axiom schema:

$$\forall y(A(y) \leftrightarrow \neg B(y)) \to A(0) \wedge \forall y < |x|\,(A(y) \to A(y+1)) \to A(|x|)$$

for $A, B \in \Sigma_1^b(L_{BA}(\mathcal{X}))$.

The following lemma is a folklore, e.g., *cf.* [17].

**Lemma 1.1** *Over $\Sigma_0^b(L_{BA}(\mathcal{X})) - LIND$ we have*

1. $\Sigma_0^b(L_{BA}(\mathcal{X})) - LCA$

2. $\Sigma_0^b(L_{BA}(\mathcal{X})) - AC$ *proves* $\Sigma_0^b(L_{BA}(\mathcal{X})) - CA, \Sigma_1^b(L_{BA}(\mathcal{X})) - AC$, $\Delta_1^b(L_{BA}(\mathcal{X})) - CA$ *and* $\Delta_1^b(L_{BA}(\mathcal{X})) - LIND$

### 1.3 The bounded arithmetic $AID$

Now we define a bounded arithmetic $AID$.

**Definition 1.1** The *language $L_{AID}$ of $AID$:* Given a linear form $\ell\|\bar{x}\|$ in $\|\bar{x}\| = \|x_1\|, \ldots, \|x_n\|$, $\Sigma_0^b$-formulae $B(\bar{x}, p), \bar{D}(\bar{x}, p) = D_1(\bar{x}, p), \ldots, D_m(\bar{x}, p)$ in $L_{BA}$ and a boolean propositional formula $I(\bar{d}, p_0, p_1)$ in atoms $\bar{d} = d_1, \ldots, d_m$ and $p_0, p_1$ we introduce an *(n+1)-ary predicate constant* $A^{\ell, B, \bar{D}, I}$. Then $L_{AID}$ is defined to be the expanded language of $L_{BA}$ having the predicate constant $A^{\ell, B, \bar{D}, I}$ for each such items $\ell, B, \bar{D}, I$. When no confusion likely occurs, we write $A$ for $A^{\ell, B, \bar{D}, I}$.

$AID$ is obtained from $\Sigma_0^b(L_{AID}) - LIND$, i.e., $\Sigma_0^b - LIND$ in the language $L_{AID}$ by adding the following *axioms for the newly introduced predicate* $A = A^{\ell, B, \bar{D}, I}$:

**(A.0)** $A(\bar{x}, p) \to 0 \neq |p| \leq \ell\|\bar{x}\|$

**(A.1)** $0 \neq |p| = \ell\|\bar{x}\| \to [A(\bar{x}, p) \leftrightarrow B(\bar{x}, p)]$

**(A.2)** $0 \neq |p| < \ell\|\bar{x}\| \to [A(\bar{x}, p) \leftrightarrow I(\bar{D}(\bar{x}, p), A(\bar{x}, p0), A(\bar{x}, p1))]$

where the LHS $I(\bar{D}(\bar{x}, p), A(\bar{x}, p0), A(\bar{x}, p1))$ denotes the result
$I(D_1(\bar{x}, p), \ldots, D_m(\bar{x}, p), A(\bar{x}, p0), A(\bar{x}, p1))$ of simultaneous replacement of the atoms $\bar{d} = d_1, \ldots, d_m, p_0, p_1$ by the formulae
$D_1(\bar{x}, p), \ldots, D_m(\bar{x}, p), A(\bar{x}, p0), A(\bar{x}, p1)$ in the boolean formula $I$, *cf.* the encoding (1).

The above **(A.0)**, **(A.1)**, **(A.2)** give the *inductive definition of the predicate* $A = A^{\ell, B, \bar{D}, I}$. To decide $A(\bar{x}, p)$ for $0 \neq |p| \leq \ell\|\bar{x}\|$ build a binary tree of depth $\ell\|\bar{x}\| - |p|$. The sons of the node $A(\bar{x}, p)$ are $A(\bar{x}, p0), A(\bar{x}, p1)$ and we put a 'gate' $I(\bar{D}(\bar{x}, p), p_0, p_1)$ there. $p$ in $A(\bar{x}, p)$ is a *clock*, i.e., it tells us the time when we stop to calculate the truth values, namely $|p| = \ell\|\bar{x}\|$ by **(A.1)**.

The point is that we can decide $A(\bar{x}, p)$ from the sons $A(\bar{x}, p0), A(\bar{x}, p1)$ propositionaly. The essence of the clause **(A.2)** is that previous stages are not quantified at all in the RHS.

We can aasign truth values $A^{\ell, B, \bar{D}, I}(\bar{x}, p)$ for $p = 0$ and for $|p| > \ell\|\bar{x}\|$ in an arbitrary manner since we need only $\{A(\bar{x}, p) : 0 \neq |p| \leq \ell\|\bar{x}\|\}$.

It is straightforward to see the following proposition, *cf.* Theorem 4 in [6].

**Proposition 1.1** *Each $\Sigma_0^b$-formula in $L_{AID}$ defines a predicate in $ALOGTIME$.*

## 2 Inductive definitions and vector summations in $AID$

In this section we show first that inductive definitions along quadtree (4-branching tree), iterated inductive definitions and simultaneous inductive definitions reduce to a single inductive definition available in $AID$. Using these we show second that vector summations are $\Sigma_0^b$-bitdefinable in $AID$.

**Lemma 2.1** *(Tree induction) For a linear form $\ell\|\bar{x}\|$ and a $\Sigma_0^b$-formula $B$ in $L_{AID}$, we have in $AID$*

$$\forall |p| \leq \ell\|\bar{x}\|[(0 \neq |p| = \ell\|\bar{x}\| \to B(p)) \,\&\, (0 \neq |p| < \ell\|\bar{x}\| \,\&\,$$
$$\bigwedge_{i<2} B(pi) \to B(p))] \to \forall |p| \leq \ell\|\bar{x}\|(0 \neq |p| \to B(p))$$

**Proof.** Apply $\Sigma_0^b - LIND$ to the following $\Sigma_0^b$ $C(u)$:

$$C(u) \Leftrightarrow_{df} \forall |p| \leq \ell\|\bar{x}\|(0 \neq |p| \,\&\, |p| + u \geq \ell\|\bar{x}\| \to B(p))$$

□

In the following lemma let $A = A^{\ell, B, \bar{D}, I}$ be a predicate defined by **(A.0)**, **(A.1)**, **(A.2)**. Let $C^+, C^-$ be $\Sigma_0^b(L_{AID})$ formulae. By separating positive, negative occurrences of atoms $p_0, p_1$ in the boolean formula $I$, we set

$$I^+(\bar{d}, p_0^+, p_0^-, p_1^+, , p_1^-) \leftrightarrow I(\bar{d}, p_0, p_1)$$

The superscript $+$ $[-]$ designates the positive [negative] occurrences. Put

$$I^-(\bar{d}, p_0^+, p_0^-, p_1^+, , p_1^-) \Leftrightarrow_{df} \neg I^+(\bar{d}, p_0^-, p_0^+, p_1^-, p_1^+)$$

Let $IH$ denote the formula:

$$\forall |p| \le \ell \|\bar{x}\| [$$
$$\{0 \ne |p| = \ell \|\bar{x}\| \to (B(\bar{x}, p) \to C^+(p)) \& (\neg B(\bar{x}, p) \to C^-(p))\} \&$$
$$\{0 \ne |p| < \ell \|\bar{x}\| \to$$
$$(I^+(\bar{D}(\bar{x}, p), C^+(p0), C^-(p0), C^+(p1), C^-(p1)) \to C^+(p)) \&$$
$$(I^-(\bar{D}(\bar{x}, p), C^+(p0), C^-(p0), C^+(p1), C^-(p1)) \to C^-(p))\}]$$

Then we see the following lemma from Lemma 2.1.

**Lemma 2.2** *(Proof by tree induction) AID proves that*

$$IH \to \forall |p| \le \ell \|\bar{x}\| [0 \ne |p| \to (A(\bar{x}, p) \to C^+(p)) \& (\neg A(\bar{x}, p) \to C^-(p))]$$

**Lemma 2.3** *(Inductive definitions along quadtrees) Let $\ell\|\bar{x}\|$ be a linear form, $B, \bar{D}$ $\Sigma_0^b$-formulae in $L_{BA}$ and $I(\bar{d}, p_{00}, p_{01}, p_{10}, p_{11})$ a boolean formulae. Define a predicate $A$ inductively by:*

**(A.0)** $A(\bar{x}, p) \to 0 \ne |p| \le 2\ell\|\bar{x}\| + 1 \& |p|$ *is odd.*

**(A.1)** *The case* $0 \ne |p| = 2\ell\|\bar{x}\| + 1$: $A(\bar{x}, p) \leftrightarrow B(\bar{x}, p)$

**(A.2)** *The case* $0 \ne |p| < 2\ell\|\bar{x}\| + 1 \& |p|$ *is odd :*

$$A(\bar{x}, p) \leftrightarrow I(\bar{D}(\bar{x}, p), \{A(\bar{x}, pij) : i, j < 2\}) \leftrightarrow$$
$$I(\bar{D}(\bar{x}, p), A(\bar{x}, p00), A(\bar{x}, p01), A(\bar{x}, p10), A(\bar{x}, p11))$$

*where* $pij = 2(2p + i) + j$

*Then $A$ can be $\Sigma_0^b$-defined in AID.*

**Proof**. For a formula $F$ and $i < 2$ put

$$F^i = \begin{cases} F & i = 1 \\ \neg F & i = 0 \end{cases} \tag{2}$$

For $i, j < 2$ put $k(ij) = Bit(2i + j, k)$.

First write the boolean formula $I$ in a DNF (disjunctive normal form):

$$\bigvee \{C_k : k < 2^4\} \leftrightarrow I(\bar{D}, \{A(\bar{x}, pij) : i, j < 2\}) \tag{3}$$

such that each disjunct $C_k$ is of the form:

$$C_k \equiv I_k(p) \wedge \bigwedge \{A(\bar{x}, pij)^{k(ij)} : i, j < 2\} \tag{4}$$

and the predicate $A$ does not occur in $I_k(p)$.

Now the 4-branching regress **(A.2)** is simulated by a tree of depth 6: first construct a $\bigvee$-tree of depth 4 corrsponding to $\bigvee$ in the DNF (3) and then construct trees of depth 2 below each leaf of $\vee$-tree to handle $\bigwedge$ in (4). We define $A(\bar{x}, p)$ using a new clock $q$ which is divided by 6-digits.

Define a predicate $A'$ inductively as follows:
**(A'.0)**
$$A'(\bar{x}, p, q) \to |p| \text{ is odd } \& \ |p| + q_0 \leq \ell + 1 \& 0 \neq |q| \leq \ell'$$
where $\ell = \ell\|\bar{x}\|$, $\ell' = 1 + 6\left(\ell - \lfloor \frac{|p|}{2} \rfloor\right)$ and

$$q_0 = \begin{cases} 2\left(\lfloor \frac{|q|\dot{-}1}{6} \rfloor\right) & \text{if } |q| \dot{-} 1 \equiv 0 \ (mod \ 6) \\ 2\left(\lfloor \frac{|q|\dot{-}1}{6} \rfloor\right) + 2 & \text{otherwise} \end{cases}$$

In the following assume $|p|$ is odd $\& \ |p| + q_0 \leq \ell + 1 \& 0 \neq |q| \leq \ell'$, the RHS of **(A'.0)**.

**(A'.1)** The case $|q| < \ell' \ \& \ |q| \equiv 1, 2, 3, 4 \ (mod \ 6)$:
$$A'(\bar{x}, p, q) \leftrightarrow A'(\bar{x}, p, q0) \vee A'(\bar{x}, p, q1)$$

**(A'.2)** The case $|q| < \ell' \ \& \ |q| \equiv 5 \ (mod \ 6)$:
$$A'(\bar{x}, p, q) \leftrightarrow \bigvee \{\exists |r| < \ell [R_2(r, q) \ \& \ I_k(p * r)] \ \& $$
$$A'(\bar{x}, p, q0) \ \& \ A'(\bar{x}, p, q1) \ \& \ k = q[0, 4) : k < 2^4\}$$

where $R_2(r, q)$ denotes the formula
$$|q| = 3(|r| \dot{-} 1) + 5 \ \& \ \forall j < 2 \forall i < \lfloor \frac{|q|}{6} \rfloor (Bit(2i + j, r) = Bit(6i + j + 4, q))$$

By Lemma 1.1.1, $\Sigma_0^b - LCA$ and Bit Extensionality axiom such a number $r$ is uniquely determined from $q$.

**(A'.3)** The case $|qi| < \ell' \ \& \ |q| \equiv 5 \ (mod \ 6)$ for an $i < 2$:
$$A'(\bar{x}, p, qi) \leftrightarrow \bigwedge \{A'(\bar{x}, p, qij)^{q(ij)} : j < 2\}$$

where $F^{q(ij)} \Leftrightarrow_{df} (q(ij) = 1 \ \& \ F) \vee (q(ij) = 0 \ \& \ \neg F)$ and $q(ij) = Bit(2i + j, q)$.

**(A'.4)** The case $|1q| = \ell' = 1 + 6\left(\ell - \lfloor \frac{|p|}{2} \rfloor\right)$:
$$A'(\bar{x}, p, 1q) \leftrightarrow \exists |r| \leq \ell [R_4(r, q) \ \& \ B(\bar{x}, p * r)]$$

where $R_4(r, q)$ denotes the formula
$$|q| = 3(|r| \dot{-} 1) \ \& \ \forall j < 2 \forall i < \lfloor \frac{|q|}{6} \rfloor (Bit(2i + j, r) = Bit(6i + j, q))$$

By induction on $6\left(\ell - \lfloor \frac{|p|}{2} \rfloor\right) - |q|$ we see the following Claim 2.1:

**Claim 2.1**

$$AID \vdash |q| \equiv 0 \ (mod\ 6)\ \&\ R_4(r,q) \to [A'(\bar{x}, p*r, 1) \leftrightarrow A'(\bar{x}, p, 1q)]$$

where $R_4(r,q)$ denotes the formula in **(A'.4)**.

By the Claim 2.1 we have for $i, j < 2$ and $1k = 2^4 + k$

$$\forall k < 2^4 (A'(\bar{x}, pij, 1) \leftrightarrow A'(\bar{x}, p, 1kij)) \tag{5}$$

Thus by putting

$$A(\bar{x}, p) \Leftrightarrow_{df} A'(\bar{x}, p, 1)$$

we get the defining axioms **(A.0)-(A.2)**. For example to see **(A.2)** assume $0 \neq |p| < 2\ell\|\bar{x}\| + 1\ \&\ |p|$ is odd. Then we have $6 < 1 + 6\left(\ell - \lfloor\frac{|p|}{2}\rfloor\right)$ by $|p|$ is odd and hence $\lfloor\frac{|p|}{2}\rfloor < \ell$. Using the defing axioms **(A'.1)-(A'.3)** of $A'$, (5), (3) and (4) we have

$$\begin{aligned}
A(\bar{x}, p) &\Leftrightarrow_{df} A'(\bar{x}, p, 1) \\
&\leftrightarrow \bigvee\{I_k(p)\ \&\ \bigwedge\{A'(\bar{x}, p, 1kij)^{k(ij)} : i, j < 2\} : k < 2^4\} \\
&\leftrightarrow \bigvee\{I_k(p)\ \&\ \bigwedge\{A'(\bar{x}, pij, 1)^{k(ij)} : i, j < 2\} : k < 2^4\} \\
&\leftrightarrow I(\bar{D}(\bar{x}, p), \{A'(\bar{x}, pij, 1) : i, j < 2\}) \\
&\Leftrightarrow I(\bar{D}(\bar{x}, p), \{A(\bar{x}, pij) : i, j < 2\})
\end{aligned}$$

$\square$

By Lemma 2.3 inductive definitions along $K$-branching trees for each constant $K \geq 2$ are $\Sigma_0^b$-definable in $AID$.

**Lemma 2.4** *(Iterated inductive definitions)*
Let $B$ and $\bar{D}$ be $\Sigma_0^b$-formulae in $L_{AID}$. Namely inductively defined predicates may occur in these formulae. Let $A = A^{\ell, B, \bar{D}, I}$ be an inductively defined predicate defined by the **(A.0)-(A.2)** in Definition 1.1 from a linear form $\ell\|\bar{x}\|$ and a boolean formula $I$. Then $A$ is $\Sigma_0^b$-definable in $AID$.

**Proof**.
**(Step1)** First we put inductively defined predicates occurring in $\bar{D}$ into the terminal condition $B$. For simplicity assume the number $m$ of atoms $d_1, \ldots, d_m$ in $I$ is 2 and let $I = I(d_2, d_3, p_0, p_1)$ and $\bar{D}(\bar{x}, p) = D_2(\bar{x}, p), D_3(\bar{x}, p)$. Further assume $\ell\|\bar{x}\| \neq 0$, i.e., $\ell\|\bar{x}\| = \sum c_i\|x_i\| + d$ with $d \neq 0$.

We say that $p$ *contains a digit* 2 *or* 3 iff $p[i, i+2) \in \{2, 3\}$ for some even $i < |p|$. Define a predicate $A'$ by induction along a quadtree as follows:
**(A'.0)** $A'(\bar{x}, p) \to |p|$ is odd & $|p| \leq 2\ell\|\bar{x}\| - 1$.
In the following assume the RHS $|p|$ is odd & $|p| \leq 2\ell\|\bar{x}\| - 1$ of **(A'.0)**.
**(A'.1)** The case $|p| < 2\ell\|\bar{x}\| - 1$ and $p$ does not contain a digit 2 nor 3:

$$A'(\bar{x}, p) \leftrightarrow I(A'(\bar{x}, p10), A'(\bar{x}, p11), A'(\bar{x}, p00), A'(\bar{x}, p01))$$

**(A'.2)** The case $|p| < 2\ell\|\bar{x}\| - 1$ and $p$ contain a digit 2 or 3:

$$A'(\bar{x}, p) \leftrightarrow A'(\bar{x}, p00)$$

**(A'.3)** The case $|p| = 2\ell\|\bar{x}\| - 1$ and $p$ does not contain a digit 2 nor 3:

$A'(\bar{x}, p)$
$\exists |q| \leq \ell\|\bar{x}\| [[|q| = \ell\|\bar{x}\| - 1 \ \& \ \forall i < |q| \dot{-} 1(Bit(i, q) = Bit(2i, p)) \ \& \ B(\bar{x}, q)]$

**(A'.4)** The case $|p| = 2\ell\|\bar{x}\| - 1$ and $p$ contains a digit 2 or 3:
**(4.1)** $A'(\bar{x}, p) \to \exists! i < |p| \ [i \text{ is even } \& \ p[i, i+2) \in \{2, 3\}]$
Let $i < |p|$ denote the unique $i$ such that $i$ is even $\& \ p[i, i+2) \in \{2, 3\}$.
**(4.k)** The case $p[i, i+2) = k \ (k \in \{2, 3\})$:

$$\begin{aligned} A'(\bar{x}, p) \quad &\leftrightarrow \quad \exists |q| \leq \ell\|\bar{x}\| [[|q| = \ell\|\bar{x}\| - 1 - \lfloor \frac{i}{2} \rfloor] \\ &\& \ \forall j < |q| \dot{-} 1 (Bit(j, q) = Bit(2j + i + 2, p)) \ \& \ D_k(\bar{x}, q)] \end{aligned}$$

Then $A'$ is defined from a boolean formula $I'$, $\Sigma_0^b$-forrmulae $\bar{D}'$ in $L_{BA}$ and a $\Sigma_0^b$-terminal condition $B'$ in $L_{AID}$ along a quadtree. From the proof of Lemma 2.3 we see that $A'$ can be defined from a boolean formula $I''$, $\Sigma_0^b$ $\bar{D}''$ in $L_{BA}$ and a $\Sigma_0^b$ $B''$ in $L_{AID}$ along a binary tree.

Further we see easily that for $k < 2$

$$A(\bar{x}, q) \leftrightarrow A'(\bar{x}, p) \text{ and } D_{1k}(\bar{x}, q) \leftrightarrow A'(\bar{x}, p1k)$$

where $10 = 2, 11 = 3$ and $p$ denotes the number such that $|p| = 2|q| - 1$ and $\forall i < |q| \dot{-} 1(Bit(2i, p) = Bit(i, q) \ \& \ Bit(2i + 1, p) = 0)$.
**(Step2)** Now we assume that no inductively defined predicate occur in $\bar{D}$. Without loss of generality we may assume that the $\Sigma_0^b$-terminal condition $B(\bar{x}, p)$ is in a prenex normal form $\forall |y| \leq \ell'\|\bar{x}\|B_0(\bar{x}, p, y)$.

Define $A'$ as follows:
**(A'.1)** The case $0 \neq |p| = \ell\|\bar{x}\| + \ell'\|\bar{x}\|$:

$$A'(\bar{x}, p) \leftrightarrow B_0(\bar{x}, p_0, p_1)$$

where $p_0 = p[\ell'\|\bar{x}\|, |p|)$ and $p_1 = p[0, \ell'\|\bar{x}\|)$.
**(A'.2)** The case $\ell\|\bar{x}\| \leq |p| < \ell\|\bar{x}\| + \ell'\|\bar{x}\|$:

$$A'(\bar{x}, p) \leftrightarrow A'(\bar{x}, p0) \wedge A'(\bar{x}, p1)$$

**(A'.3)** The case $0 \neq |p| < \ell\|\bar{x}\|$:

$$A'(\bar{x}, p) \leftrightarrow I(\bar{D}(\bar{x}, p), A'(\bar{x}, p0), A'(\bar{x}, p1))$$

Then $|p| = \ell\|\bar{x}\| \to (A'(\bar{x}, p) \leftrightarrow \forall |y| \leq \ell'\|\bar{x}\|B_0(\bar{x}, p, y))$ and hence

$$0 \neq |p| \leq \ell\|\bar{x}\| \to (A(\bar{x}, p) \leftrightarrow A'(\bar{x}, p))$$

The number of sharply bounded quantifiers in the terminal condition $B_0$ for $A'$ is fewer than one in $B$ for $A$. Thus we may assume that the terminal condition $B(\bar{x}, p)$ is an open formula in DNF. By considering a little bit higher tree, *cf.* proof of Lemma 2.3, we may assume further that the terminal condition is a literal.

(**Step3**) Now suppose that the terminal condition $B(\bar{x}, p)$ is an atomic formula of the form $A'(\bar{t}(\bar{x}, p), s(\bar{x}, p))$ (or its negation) for some terms $\bar{t}, s$ and an inductively defined predicate $A'$. Namely $A$ is defined by:

(**A.2**) $0 \neq |p| < \ell\|\bar{x}\| \to [A(\bar{x}, p) \leftrightarrow I(\bar{D}(\bar{x}, p), A(\bar{x}, p0), A(\bar{x}, p1))]$

(**A.1**) $0 \neq |p| = \ell\|\bar{x}\| \to [A(\bar{x}, p) \leftrightarrow A'(\bar{t}(\bar{x}, p), s(\bar{x}, p))]$

While $A'$ is defined from $\ell', I'$ and some $\Sigma_0^b$-formulae $B', \bar{D}'$ in $L_{BA}$ as follows:

(**A'.1**) $0 \neq |p| = \ell'\|\bar{y}\| \to [A'(\bar{y}, p) \leftrightarrow B'(\bar{y}, p)]$

(**A'.2**) $0 \neq |p| < \ell'\|\bar{y}\| \to [A'(\bar{y}, p) \leftrightarrow I'(\bar{D}'(\bar{y}, p), A'(\bar{y}, p0), A'(\bar{y}, p1))]$

Define $A$ alternatively as follows:

(**A.2**) $0 \neq |p| < \ell\|\bar{x}\| \to [A(\bar{x}, p) \leftrightarrow I(\bar{D}(\bar{x}, p), A(\bar{x}, p0), A(\bar{x}, p1))]$.
  In the following assume $0 \neq |p| \geq \ell\|\bar{x}\|$. Put $p = q * r$ with $|q| = \ell\|\bar{x}\|$ and $r > 0$. Our aim is to have

$$A(\bar{x}, p) \leftrightarrow A'(\bar{t}(\bar{x}, q), s(\bar{x}, q) * r) \qquad (6)$$

(**a**) $A(\bar{x}, p) \to 0 \neq |s(\bar{x}, q)| + |r| \dot{-} 1 \leq \ell'\|\bar{t}(\bar{x}, q)\|$.
  Assume the RHS $0 \neq |s(\bar{x}, q)| + |r| \dot{-} 1 \leq \ell'\|\bar{t}(\bar{x}, q)\|$ of (**a**).

(**b**) The case $|s(\bar{x}, q)| + |r| \dot{-} 1 < \ell'\|\bar{t}(\bar{x}, q)\|$:

$$A(\bar{x}, p) \leftrightarrow I'(\bar{D}'(\bar{t}(\bar{x}, q), s(\bar{x}, q) * r), A(\bar{x}, p0), A(\bar{x}, p1)$$

(**c**) The case $|s(\bar{x}, q)| + |r| \dot{-} 1 = \ell'\|\bar{t}(\bar{x}, q)\|$:

$$A(\bar{x}, p) \leftrightarrow B'(\bar{t}(\bar{x}, q), s(\bar{x}, q) * r)$$

By these (**a**),(**b**),(**c**) we see (6). Also by $|q| = \ell\|\bar{x}\|$, $|p|$ varies through

$$|p| \leq \ell\|\bar{x}\| + \ell'\|\bar{t}(\bar{x}, q)\| - |s(\bar{x}, q)| \leq \ell\|\bar{x}\| + \ell''\|\bar{x}\|$$

for some linear form $\ell''$. □

Next we show that simualtaneous inductive definitions and vector summations are definable in $AID$.

**Lemma 2.5** *(Simultaneous induction)*
*Let $\ell = \ell\|\bar{x}\|$ be a linear form, $B(\bar{x}, \lambda, p), \bar{D}(\bar{x}, \lambda, p)$ $\Sigma_0^b$-formulae in $L_{BA}$ and $I(\bar{d}, \{p_{ji} : j \leq K, i < 2\})$ be a boolean formula with a constant $K$. Define a predicate $A$ inductively by:*

**(A.0)** $A(\bar{x}, \lambda, p) \to 0 \neq |p| \leq \ell \,\&\, |\lambda| \leq \ell$.
 Assume the RHS $0 \neq |p| \leq \ell \,\&\, |\lambda| \leq \ell$ in the following.

**(A.1)** $|p| = \ell \to [A(\bar{x}, \lambda, p) \leftrightarrow B(\bar{x}, \lambda, p)]$

**(A.2)** $|p| < \ell \to [A(\bar{x}, \lambda, p) \leftrightarrow I(\bar{D}(\bar{x}, \lambda, p), \{A(\bar{x}, \lambda + j, pi) : i < 2, j \leq K\})$

Then $A$ is $\Sigma_0^b$-definable in $AID$.

**Definition 2.1** (Bitdefinable functions)
Let $T$ be a sound theory and $\Phi$ be a set of formulae in the language of $T$. We say that a function $f(\bar{x})$ of polynomial growth rate is $\Phi$-*bitdefinable* if its *bitgraph* $\{i < |f(\bar{x})| : Bit(i, f(\bar{x})) = 1\}$ is definable by a formula in $\Phi$, i.e., there exists a formula $A_f(i, \bar{x})$ in $\Phi$ such that in the standard model $A_f(i, \bar{x}) \leftrightarrow Bit(i, f(\bar{x})) = 1$. If the theory $T$ proves the axiom schema $\Phi - CA$ and the Bit Extensionality axiom, then $T \vdash \forall \bar{x} \exists ! y (y = \{i < p(|\bar{x}|) : A_f(i, \bar{x})\})$ for a polynomial $p$ with $|f(\bar{x})| \leq p(|\bar{x}|)$.

**Theorem 2.1** (Bounded vector summation)
If $f(i, \bar{y})$ is $\Sigma_0^b$-bitdefinable in $AID$, then so is the function

$$g(x, \bar{y}) = \sum_{i < |x|} f(i, \bar{y})$$

Also the defined function $g$ enjoys demonstrably in $AID$,
$g(0, \bar{y}) = 0$ and $g(x, \bar{y}) = g(\lfloor \frac{x}{2} \rfloor, \bar{y}) + f(|x| - 1, \bar{y})$ for $x \neq 0$.

**Corollary 2.1** (Bounded counting)
For each $\Sigma_0^b$-formula $\varphi$ in $L_{AID}$, the function

$$C_\varphi(x) = \#\{i < |x| : \varphi(i)\}$$

is $\Sigma_0^b$-bitdefinable in $AID$. We always have $C_\varphi(x) \leq |x|$. Therefore we can use freely the function $C_\varphi$ in $\Sigma_0^b$-formulae.

First assuming Lemma 2.5 we show Theorem 2.1 using carry-save-addition to combine four numbers into two, in p.922 of [5].

**Proof** of Theorem 2.1. For simplicity we suppress parameters $\bar{y}$. Pick a polynomial $H$ so that $|\sum_{i<|x|} f(i)| \leq H(|x|)$.

Define a predicate $A(x, \lambda, p)$ by simultaneous recursion in Lemma 2.5 as follows:

**(A.0)** $A(x, \lambda, p) \to 0 \neq |p| \leq \|x\| + 1 \,\&\, p < 2^{\|x\|} + |x| \,\&\, \lambda \leq 2H(|x|)$.
 Assume the RHS of **(A.0)**. In what follows we write

$$\begin{aligned} s_\lambda^p &\Leftrightarrow_{df} S(x, \lambda, p) \Leftrightarrow_{df} A(x, 2\lambda, p) \\ c_\lambda^p &\Leftrightarrow_{df} C(x, \lambda, p) \Leftrightarrow_{df} A(x, 2\lambda + 1, p) \end{aligned}$$

(**A.1**) The case $|p| = \|x\| + 1$:

$$\begin{aligned} s_\lambda^p &\leftrightarrow A_f(\lambda', p') (\Leftrightarrow Bit(\lambda', f(p')) = 1) \\ c_\lambda^p &\leftrightarrow \bot \end{aligned}$$

where $\lambda + \lambda' = H(|x|)$ and $p = 2^{\|x\|} + p', p' < |x|$, i.e., reverse the digits in $f(p')$.

(**A.2**) The case $|p| < \|x\| + 1$:

$$\begin{aligned} s_\lambda^p &\leftrightarrow (s_\lambda^{p0} \oplus s_\lambda^{p1} \oplus c_\lambda^{p0}) \oplus (\#\{s_{\lambda+1}^{p0}, s_{\lambda+1}^{p1}, c_{\lambda+1}^{p0}\} \geq 2) \oplus c_\lambda^{p1} \\ c_\lambda^p &\leftrightarrow \#\{s_{\lambda+1}^{p0}, s_{\lambda+1}^{p1}, c_{\lambda+1}^{p0}, \#\{s_{\lambda+2}^{p0}, s_{\lambda+2}^{p1}, c_{\lambda+2}^{p1}\} \geq 2, c_{\lambda+1}^{p1}\} \geq 2 \end{aligned}$$

where $\oplus$ denotes the excluded or, and for propositions $p, q, r$,

$$\#\{p, q, r\} \geq 2 \Leftrightarrow_{df} (p \wedge q) \vee (q \wedge r) \vee (r \wedge p)$$

For each fixed $p = 2^{\|x\|} + p', p' < |x|$ $s^p = \{s_\lambda^p\}_{\lambda < H(|x|)}$ gives the reverse binary representation of a number, and similarly for $c^p = \{c_\lambda^p\}_{\lambda < H(|x|)}$.

(**A.1**) $s^p$ is the reverse binary representation of the number $f(p')$.

(**A.2**) In reverse notation, first by carry-save-addition to combine three numbers into two, $s^{p0} + s^{p1} + c^{p0} = s' + c'$, and then $s' + c' + c^{p1} = s^p + c^p$ once again by carry-save-addition. Thus $s^{p0} + s^{p1} + c^{p0} + c^{p1} = s^p + c^p$.

Let $s_\lambda, c_\lambda$ denote the following formulae with $\lambda + \lambda' = H(|x|)$

$$s_\lambda \leftrightarrow s_{\lambda'}^1 \; ; \; c_\lambda \leftrightarrow c_{\lambda'}^1$$

Let $a_\lambda$ denote the full addition of $s_\lambda$ and $c_\lambda$:

$$a_\lambda \leftrightarrow s_\lambda \oplus c_\lambda \oplus \exists \eta < \lambda (s_\eta \wedge c_\eta \wedge \forall \rho \in (\eta, \lambda)(s_\rho \oplus c_\rho))$$

Then $a_\lambda$ is a $\Sigma_0^b$-formula such that

$$a_\lambda \leftrightarrow Bit(\lambda, \sum_{p' < |x|} f(p')) = 1$$

for $\lambda \leq H(|x|)$. $\square$

**Proof** of Lemma 2.5. We show that $A(\bar{x}, \lambda, 1)$ is $\Sigma_0^b$-definable. We define a predicate $A'$ so that, for $|p| \equiv 1 \, (mod \, 1 + K)$,

$$A'(\bar{x}, \lambda, p) \leftrightarrow A(\bar{x}, \lambda + \#(J), q)$$

where $|p|= 1 + r(1+K)$ and $|q|= 1 + r, |J|\leq rK$, and their digits are defined by $Bit(i,q) = Bit(K + i(1+K), p)$ for $i < r$, and $Bit(i, J) = Bit(i + \lfloor \frac{i}{K} \rfloor, p)$ for $i < rK$.

$\#(J)$ denotes the number

$$\#(J) = \#\{i < |J|: Bit(i, J) = 1\}$$

Then

$$A(\bar{x}, \lambda, 1) \leftrightarrow A'(\bar{x}, \lambda, 1)$$

Such an $A'$ is defined as follows: Each $j \leq K$ is coded by $0^{[K-j]}1^{[j]}$ in unary notation. Put $\ell = \ell \|\bar{x}\|$ and $\ell' = 1 + (\ell - 1)(1+K)$.

**(A'.0)** $A'(\bar{x}, \lambda, p) \to 0 \neq |p| \leq \ell' \ \& \ |p| \equiv 1 \ (mod \ 1 + K) \ \& \ |\lambda| \leq \ell$.

Assume the RHS of **(A'.0)** in the following. Also let $q$ and $J$ denote the numbers defined above.

**(A'.1)** The case $|p|= \ell'$:

$$A'(\bar{x}, \lambda, p) \leftrightarrow B(\bar{x}, \lambda + \#(J), q)$$

**(A'.2)** The case $|p|< \ell'$:

$$A'(\bar{x}, \lambda, p) \leftrightarrow I(\bar{D}(\bar{x}, \lambda, q), \{A'(\bar{x}, \lambda, pi0^{[K-j]}1^{[j]}) : i < 2, j \leq K\})$$

In the RHS $A'(\bar{x}, \lambda, pi0^{[K-j]}1^{[j]})$ is substituted in place of $A(\bar{x}, \lambda + j, pi)$ in **(A.2)**. The definition is along a $2(1+K)$-branching tree.

Thus the problem reduces to show the

**Proposition 2.1** $\#(J)$ is $\Sigma_0^b$-bitdefinable for $|J|\leq \|\bar{x}\|$

This is a bounded counting with a bound $\|\bar{x}\|$, while Corollary 2.1 is a bounded counting with a bound $|x|$. By repeating the above proofs Proposition 2.1 reduces to show the

**Proposition 2.2** $\#(J)$ is $\Sigma_0^b$-definable for $|J|\leq c |x|_3 +c$ for any constant $c$.

**Proof** of Proposition 2.2. Suppose $|J|\leq c |x|_3 +c$. It suffices to show that $y \leq \#(J)$ is $\Sigma_0^b$-definable. Define

$$y \leq \#(J) \ \Leftrightarrow_{df} \ \exists |u|\leq \|J\| \cdot y + 1 \forall i < y[Bit(u_i, J) = 1 \ \& $$
$$(i + 1 < y \to u_i < u_{i+1}) \ \& \ y \leq |J|]$$

where $u_i = u[\|J\| \cdot i, \|J\| \cdot (i + 1))$. This means that for some $u_i$, $u = u'_{y-1} * \cdots * u'_1 * u'_0$ for $u'_i = 2^{\|J\|} + u_i$ and $u_0 < u_1 < \cdots < u_{y-1} < |J| \ \& \ \{u_i : i < y\} \subseteq \{v < |J|: Bit(v, J) = 1\}$. Each $u_i <|J|$ and hence $\|J\|$-digits suffices to represent $u_i$ in binary notation. Also note that any multiplication occurring in this definition is a multiplication for small numbers, cf. $multi(i, j, x, y)$ in section 1.

The quantifier $\exists |u|\leq \|J\| \cdot y$ with $y \leq |J|$ is a sharply bounded one since $\|J\|\cdot |J|\leq \ell\|x\|$ for some linear form $\ell$ by the supposition $|J|\leq c |x|_3 +c$. □

# 3  $ALOGTIME$ are $\Sigma_0^b$-definable in $AID$

In this section we show that each predicate in $ALOGTIME$ is $\Sigma_0^b$-definable in $AID$.

**Theorem 3.1** *Each predicate in $ALOGTIME$ is $\Sigma_0^b$-definable in $AID$.*

**Proof**. (*cf*. p.64, p.77 in [2] for indexing alternating Turing machines.) Let $A$ be a predicate in $ALOGTIME$ and $M = (Q, \Sigma, \delta, q_0, g)$ be an alternating Turing machine which recognizes $A$ such that $M$ always halts in time $\ell\|x\|$ on input $x$ for a linear form $\ell$. We assume that

1. $Q$ is a finite set of states and $q_0 \in Q$ is an initial state.

2. $\Sigma = \{0,1\}, \bar{\Sigma} = \{0,1,b\}$, where $b$ denotes the blank.

3. $g : Q \to \{\wedge, \vee, accept, reject\}$.

4. $\delta$ is a transition function such that $\delta : Q \times \bar{\Sigma}^{k+2} \to \mathcal{P}(\bar{\Sigma}^{k+1} \times H^{k+1} \times Q)$ with $H = \{L, N, R\}$.

5. The meanings of $L, R, N$ are given as follows. $L$: moving one cell to the left, $R$: moving one cell to the right, $N$: do not move.

6. $M$ contains $k+2$-tapes; a read-only input tape, $k$-work tapes and an index tape. $M$ writes down a number $i \leq |x|$ for an input $x$ in binary notation on the index tape to read the $i^{th}$ input symbol on the input tape. These tapes are numbered in this order. Thus the input tape is referred as the 0-th tape.

Without loss of generality we can assume that the computation tree of $M$ on input $x$ is a binary tree of depth $\ell\|x\|$. Each $w \in \Sigma^*$ corresponds to a node in a computation tree. For a $w \in \Sigma^*$ let $pd(w) \in \Sigma^*$ denote a word such that $w = pd(w)j$ for some $j \in \Sigma$, i.e., $pd(w)$ is obtained from $w$ by deleting the rightmost symbol in $\Sigma$. Put $w_0 \subset w_1 \Leftrightarrow_{df} w_0$ is an initial segment of $w_1$ for words $w_0, w_1$.

Thus we have functions (definable by some terms in the language $L_{BA}$) $\delta^\Sigma, \delta^H, \delta^Q$ such that for $q \in Q, \bar{s} \in \bar{\Sigma}^{k+2}, 1 \leq j \leq k+1, w \in \Sigma^*$ with $|w| \leq \ell\|x\|$

$$(\delta^\Sigma(q, \bar{s}, j, w) : 1 \leq j \leq k+1) * (\delta^H(q, \bar{s}, j, w) : 1 \leq j \leq k+1) * (\delta^Q(q, \bar{s}, w))$$

is in $\delta(q, \bar{s})$.

This $2k+3$-tuple denotes the next move at $w$ when at the predecessor node $pd(w)$, the state is $q$ and the scanned symbols are $\bar{s}$. $\delta^\Sigma(q, \bar{s}, k+1, w)$ is the symbol written on the index tape.

Let $\bar{q}, \bar{H}, \bar{s}$ denote the following objects:

1. $\bar{q} = (q^i : i \leq \ell\|x\|) \in Q^*, q^i \in Q$

2. $\bar{H} = H_1, \ldots, H_{k+1}$ and for each $j$ with $1 \leq j \leq k+1$,
   $H_j = (H_j^i : i \leq \ell\|x\|) \,\&\, H_j^i \in \{L, N, R\}$.

3. $\bar{s} = s_0, s_1, \ldots, s_{k+1}$ and for each $j \leq k+1$,
   $s_j = (s_j^i : i \leq \ell\|x\|) \,\&\, s_j^i \in \{0, 1, b\}$.

Let $I \in \{0,1\}^*$ denote a path through the computation tree, $|I| = \ell\|x\|$. Then these objects denote guesses on $I$:

1. $q^i$ is a guess of the state at node $w \subset I$ with $|w| = j$.

2. $H_j$ is a guess of the moves of the $j^{th}$ head on $I$.

3. $s_j$ is a guess of the scanned symbols by the $j^{th}$ head on $I$.

For $j$ with $1 \leq j \leq k+1$, $|w| \leq \ell\|x\|$ and $\lambda \leq |x|$, $C_j(x, \bar{q}, \bar{H}, \bar{s}, \lambda, w)$ denotes the symbol in $\bar{\Sigma}$ written on the $\lambda^{th}$ cell of $j^{th}$ tape at node $w$ when $\bar{q}, \bar{H}, \bar{s}$ are guesses on the path $I$: First for the empty word $\varepsilon$, $C_j(x, \bar{q}, \bar{H}, \bar{s}, \lambda, \varepsilon) = b$. Next suppose $w \neq \varepsilon$. Let $Position(H_j, w)$ denote the position of $j^{th}$ head at node $w$. $Position(H_j, w)$ is determined by counting the numbers of $L$'s and $R$'s in the first $|w|$ part of $H_j$. Thus $Position(H_j, w)$ is $\Sigma_0^b$-definable.

**Case1** $\forall w_1 \subset w(\lambda \neq Position(H_j, w_1))$: $\lambda$ is not and has not been the position where $j^{th}$ head stays at node $w$ or has stayed before $w$. Then $C_j(x, \bar{q}, \bar{H}, \bar{s}, \lambda, i) = b$.

**Case2** $\exists w_1 \subset w(\lambda = Position(H_j, w_1))$: Let $w_1$ denote the latest such node. Letting $w_0 = pd(w_1)$, $C_j(x, \bar{q}, \bar{H}, \bar{s}, \lambda, w) = \delta^\Sigma(q^{|w_0|}, \bar{s}^{|w_0|}, j, w_1)$ with $\bar{s}^{|w_0|} = (s_j^{|w_0|} : 0 \leq j \leq k+1)$.

Thus $C_j(x, \bar{q}, \bar{H}, \bar{s}, \lambda, i)$ is $\Sigma_0^b$-definable in $L_{AID}$.

Letting $w_0 = pd(w)$ put

$$Init(\bar{q}, \bar{H}, \bar{s}) \Leftrightarrow_{df} q^0 = q_0 \text{ (initial state)} \,\&\, \bigwedge\{s_j^0 = b : j \leq k+1\}$$

$$State(\bar{q}, \bar{s}, I) \Leftrightarrow_{df} \forall w \subset I[w \neq \varepsilon \Rightarrow q^{|w|} = \delta^Q(q^{|w_0|}, \bar{s}^{|w_0|}, w)]$$

$$Head(\bar{q}, \bar{H}, \bar{s}, I) \Leftrightarrow_{df} \forall w \subset I[w \neq \varepsilon \Rightarrow \bigwedge_{1 \leq j \leq k+1} H_j^{|w|} = \delta^H(q^{|w_0|}, \bar{s}^{|w_0|}, j, w)]$$

$$Symbol(\bar{q}, \bar{H}, \bar{s}, I) \Leftrightarrow_{df} \forall w \subset I[w \neq \varepsilon \Rightarrow s_0^{|w|} = Bit(C_{k+1}^w, x) \,\&\,$$
$$\bigwedge_{1 \leq j \leq k+1} s_j^{|w|} = C_j(x, \bar{q}, \bar{H}, \bar{s}, \lambda, w)]$$

where $\lambda = Position(H_j, w)$ and $C_{k+1}^w$ denotes a number in binary notation, i.e., $C_{k+1}^w \in \{0,1\}^*$ such that $Bit(\lambda, C_{k+1}^w) = 1 \Leftrightarrow C_{k+1}(x, \bar{q}, \bar{H}, \bar{s}, \lambda, w) = 1$, in other words $C_{k+1}^w$ is the content of the $(k+1)^{th}$ index tape at the node $w$.

Let $q(w) = q(w, I)$ denote the state at a node $w \subset I$:

$$q(w) = q \Leftrightarrow \exists \bar{q} \exists \bar{H} \exists \bar{s}[Init(\bar{q}, \bar{H}, \bar{s}) \,\&\, State(\bar{q}, \bar{s}, I) \,\&\, Head(\bar{q}, \bar{H}, \bar{s}, I)$$
$$\&\, Symbol(\bar{q}, \bar{H}, \bar{s}, I) \,\&\, q^{|w|} = q]$$

$\bar{q}, \bar{H}, \bar{s}$ are words on length at most $\ell' \|x\|$ for a linear $\ell'$ over a finite alphabets $Q, H = \{L, N, R\}, \bar{\Sigma} = \{0, 1, b\}$, resp. Therefore these existential quantifiers are sharply bounded and hence $q(w)$ is $\Sigma_0^b$-definable. Further existence and uniqueness conditions for $\bar{q}, \bar{H}, \bar{s}$ are provable from $\Sigma_0^b - LIND$.

Define a $\Sigma_0^b$-predicate $A_M(x, w)$ in $L_{AID}$ such that

$$|w| = \ell \|x\| \to [A_M(x, w) \leftrightarrow q(w) \text{ is an accepting state, i.e., } g(q(w)) = accept$$

and

$|w| < \ell \|x\| \to [A_M(x, w) \leftrightarrow$

$[q(w)$ is a universal state, i.e., $g(q(w)) = \wedge \,\&\, A_M(x, w0) \,\&\, A_M(x, w1)] \vee$

$[q(w)$ is an existential state, i.e., $g(q(w)) = \vee \,\&\, (A_M(x, w0) \vee A_M(x, w1))]$

Thus the given predicate $A \in ALOGTIME$ is defined by $A(x) \leftrightarrow A_M(x, \varepsilon)$. □

**Remark.** If we donot guess $\bar{q}, \bar{H}, \bar{s}$ and define directly the configuration $C(x, i)$ on input $x$ at the node $i$, then the resulting definition would involve a complicated simultaneous inductive definition, even if $M$ is a deterministic Turing machine with run times at most $\ell \|x\|$.

## 4 Consistency proof of Frege system

In this section we show that a truth definition for PLOF formulae (Postfix-Longer-Operands-First) is $\Sigma_0^b$-definable in $AID$. This is done by mimicing the proofs in Buss [7] almost word for word. The reader is recommended to have a copy of [7] in hand. Although one could apply the simplified algorithm for boolean formula evaluation in [8], we stick to [7], since the latter gave a full proof of the fact that the truth definition respects the meaings of propositional proofs.

In the next section we show that, if $f(\bar{x})$ is a $\Sigma_0^b$-bitdefinable function in $AID$, then $C(f(\bar{x}))$ is $\Sigma_0^b$-definable in $AID$ for any $\Sigma_0^b$-formula $C(y)$, cf. Definitions 5.5, 5.6 and Lemma 5.10. Therefore we can use freely such functions in $\Sigma_0^b$-formulae.

Let $x$ be (a code of) a sequence of $19 < 2^5$ symbols in

$$\Sigma = \{p, 0, 1, (,)(\text{parentheses}), , (\text{comma}),$$
$$13 \text{ propositional connectives (unary or binary)}\}$$

cf. p.8 [7].:

1. the length $|x|_\Sigma$ of $x$ as a word from $\Sigma$ 5 $|x|_\Sigma = |x|$,

2. $j^{th}$ symbol from $\Sigma$ in $x$ $Sym_j^x = x[5(j-1), 5j)$ for $1 \leq j \leq |x|_\Sigma$.

3. A *logical symbol* is a parenthesis, comma, propositional connective or propositional variable $p_i$ ($i$ is a word on $\{0, 1\}$),

4. $x[i]$ = the $i^{th}$ logical symbol of $x = x[5m, 5n)$, where $m = \min\{j \leq |x|: Sym_j^x$ is not in $x[i]\} + 1$ and $n = \max\{j \leq |x|: Sym_j^x$ is in $x[i+1]\}$, where

5. $Sym_j^x$ is in $x[i] \Leftrightarrow_{df} i = \#\{k \leq j : Sym_k^x$ is not 0 nor 1$\}$.

Note that $\min\{j \leq |x|: \cdots\}, \max\{j \leq |x|: \cdots\}$ with $\Sigma_0^b$ conditions $\cdots$ are available in $\Sigma_0^b$-formulae by $\Sigma_0^b - LIND$.

$|x|_L$ = the number of logical symbols in $x = \#\{j \leq |x|: Sym_j^x$ is not 0 nor 1$\}$

First of all we have to develop metamathematics, arithmetization of syntax, e.g., define

$x$ is (a code of) a postfix formula $\Leftrightarrow x[1]$ is an atomic formula &

$\#i \leq |x|_L$ ($x[i]$ is an atomic formula) =

$1 + \#i \leq |x|_L$ ($x[i]$ is a binary connective)

$\forall j < |x|_L$ [$\#i \leq j(x[i]$ is an atomic formula$) >$

$\#i \leq j(x[i]$ is a binary connective$)$]

This requires *countings*, $C_\varphi(x) = \#\{i < |x|: \varphi(i)\}$.

It is easy to see the

**Proposition 4.1** *There exists a $\Sigma_0^b$-formula $PL(x, i)$ such that if $x$ is an infix formula, then $y = \{i < |x|: PL(x, i)\}$ is the PLOF form of $x$, i.e., for any $k < |x|_L$, the $k^{th}$ symbol of $y$ is the $k^{th}$ PLOF symbol of $x$, cf. p.12 [7].*

In the following, otherwise stated, $x$ denotes a PLOF formula. By $\Sigma_0^b - LIND$, which corresponds to brute forth induction in [7], we have, cf. p.16 [7]: $\forall j \leq |x|_L \exists! i \leq |x|_L \{x[i, j]$ is a formula$\}$, where $x[i, j]$ denotes a substring of $x$ from $x[i]$ through $x[j]$ inclusive. Therefore we have a $\Sigma_0^b$-definable function $x_j$ such that

$x_j$ = the unique subformula of $x$ of the form $x[i, j]$

For $i, j \leq |x|_L$, cf. p.16 [7],

$j \trianglelefteq i \Leftrightarrow_{df} x[j]$ is in $x_i \Leftrightarrow_{df} k \leq j \leq i$ with $x_i = x[k, i]$

$lca(j, i) =_{df} \min\{k \leq |x|_L : i \trianglelefteq k \,\&\, j \trianglelefteq k\}$

Let $x[i, j]$ be a $\leq 1$-scarred formula, cf. p.15 [7]. Suppose $l < r, l < j$ and $i \leq r$. Then $k$ is the *breakpoint* of $x[i, j]$ 1-*selected by* $(l, r)$ if

$k = \max\{k \leq \min\{r, j\} : \text{either } x[l+1] \text{ or } x[i] \text{ is in } x_k\}$, cf. p.16 [7]. This is a $\Sigma_0^b$-definition.

Define, cf. p.16 [7], $\Delta_u$ and $\varepsilon_u$ inductively: $\Delta_0 = 2, \varepsilon_u \lfloor \frac{1}{2} \Delta_u \rfloor, \Delta_{u+1} = \Delta_u + \varepsilon_u$, i.e., $\Delta_u = \lfloor \frac{3}{2} \lfloor \frac{3}{2} \cdots \lfloor \frac{3}{2} 2 \rfloor \cdots \rfloor \rfloor$ with $u$'s $\lfloor \frac{3}{2} \cdot \rfloor$.

$\Delta_u$ is needed to be defined up to $|x|_L < \Delta_{u-1}$. Therefore $u < \|x\|$ suffices.

To see that $\Delta_u$ $(u < \|x\|)$ is $\Sigma_0^b$-definable, we use the carry-save-addition as in the proof of Theorem 2.1: for $|p| \leq \|x\|$, the node $p$ codes the number $\{S(x, \lambda, p) + C(x, \lambda, p)\}_\lambda$, and the number corresponds to $\Delta_{\|x\|-|p|}$ if $p$ is even, and to $\varepsilon_{\|x\|-|p|}$ if $p$ is odd. For $u < \|x\|$, $\Delta_u < (\frac{3}{2})^{u+2} < 2^{u+2} \leq 8|x|$. Thus we can use the functions $\Delta_u$ and $\varepsilon_u$ in $\Sigma_0^b$-formulae.

One can $\Sigma_0^b$-define the followings in pp.17-19 [7]:

1. Breakpoints $a_p$ $(p = 1, 2, 3, 4)$ of $x[i, j]$ generated by $(m, n)$ with $n - m = \Delta_{u+1}$ for some $u \geq 0$.

   (a) $a_1$ is the breakponit of $x[i, j]$ 1-selected by $(m, m + \varepsilon_u]$
   
   (b) $a_2$ is the breakponit of $x[i, j]$ 1-selected by $(m + \varepsilon_u, n - \varepsilon_u]$
   
   (c) $a_4$ is the least common ancestor $lca(a_1, a_2)$ of $a_1, a_2$
   
   (d) $a_3 = a_4 - 1$

2. For a formula $x[i, j]$ and numbers $n, m$ such that $m < i \leq j \leq n$ and $n - m = \Delta_{u+1}$, split the formula $x[i, j]$ into up to $\leq 4$ subformulae $SubFm_1(x, [i, j], (m, n), p)$ $(p = 1, 2, 3, 4)$ by introducing the breakpoints $a_p$ generated by $(m, n)$:

$$SubFm_1(x, [i, j], (m, n), 1)] = [i, a_1]$$
$$SubFm_1(x, [i, j], (m, n), 2) = [a_1 + 1, a_2]$$
$$SubFm_1(x, [i, j], (m, n), 3) = [a_2 + 1, a_3]$$
$$SubFm_1(x, [i, j], (m, n), 4) = [a_4 + 1, j]$$

   where breakpoints are defined so that, if $a_2 \neq a_4$, then $BinOp(x, [i, j], (m, n)) = x[a_4]$ is a binary connective.

3. $x[i]$ is a *scar* of the interval $[a, b]$ iff $i < a$ and there is a connective $x[k]$ with $a \leq k \leq b$ such that $x_i$ is one of the operands of $x[k]$.

Then Lemma 12 in p.19 [7] is provable in $AID$.

**Lemma 4.1** (=*Lemma 12 in p.19 [7]*) *AID proves the followings: let $x[i, j]$ be a $\leq 1$-scarred formula, $n - m = \Delta_{u+1}$ and $a_p$ $(p = 1, 2, 3, 4)$ the breakpoints of $x[i, j]$ generated by $(m, n)$.*

**(a)** $i \leq m + \varepsilon_u < j \Rightarrow a_1 + 1 \trianglelefteq a_2$

**(b)** $\max\{m+1, i\} \leq a \leq a_2 \Rightarrow a \trianglelefteq a_1 \vee a \trianglelefteq a_2$ *and* $\max\{m+1, i\} \leq a \leq a_3 \Rightarrow a \trianglelefteq a_1 \vee a \trianglelefteq a_2 \vee a \trianglelefteq a_3$

(c) *For $p = 1, 2, 3, 4$, $SubFm_1(x, [i,j], (m,n), p)$ does not have more than one scar $x[k]$ with $k \geq \max\{m+1, i\}$*

In the proof of Lemma 4.1.(a) use the fact that $x$ is a PLOF formula. For a proof of Lemma 4.1 we need the

**Proposition 4.2**  *1. $a \trianglelefteq b \Rightarrow a \leq b$*

*2. $x_b = x[a,b] \Rightarrow (a \leq c \leq b \Leftrightarrow c \trianglelefteq b)$*

*3. $a < b \ \& \ a \perp b (\Leftrightarrow_{df} a \not\trianglelefteq b \ \& \ b \not\trianglelefteq a) \ \& \ c = lca(a,b) \Rightarrow b \trianglelefteq c - 1$.*

*cf.* p.21 [7]. Let $n - m = \Delta_{u+1}$. For $k \leq u+1$ and $1 \leq p_1, \ldots, p_k \leq 4$ ($p_1, \ldots, p_k$ is coded by a number of length $1 + 2k \leq 3 + 2u < 3 + 2\|x\|$) define

$$Int_k((m,n], p_1, \ldots, p_k) = (m', n']$$

with

$$m' = m + \sum_{j=1}^{k} \lfloor \frac{1}{2}(p_j - 1) \rfloor \cdot \varepsilon_{u+1-j} \text{ and } n' = m' + \Delta_{u+1-k}$$

and

$$\lfloor \frac{1}{2}(p - 1) \rfloor = \begin{cases} 0 & \text{if } p = 1, 2 \\ 1 & \text{if } p = 3, 4 \end{cases}$$

This can be $\Sigma_0^b$-defined by using *vector summation* $g(x, \bar{y}) = \sum_{i < |x|} f(i, \bar{y})$.

**Definition 4.1** Iteration of splitting into subformulae $SubFm_k$. *cf.* p.22 [7]. For $k \leq u+1, 1 \leq p_1, \ldots, p_k \leq 4, 1 \leq l \leq k$ let $a_1^l, \ldots, a_4^l$ be breakponits of $x[i,j]$ generated by $Int_{l-1}((m,n], p_1, \ldots, p_{l-1})$. Put $a_0^l = i - 1, a_5^l = j$ and

$$i_l = \begin{cases} p_l - 1 & \text{if } 1 \leq p_l \leq 3 \\ 4 & \text{if } p_l = 4 \end{cases}$$

Then $SubFm_1(x, [i,j], Int_{l-1}((m,n], p_1, \ldots, p_{l-1}), p_l) = [a_{i_l}^l + 1, a_{1+i_l}^l]$. Put

$$c_k = \max\{a_{i_l}^l : 1 \leq l \leq k\}, \ d_k = \min\{a_{1+i_l}^l : 1 \leq l \leq k\}$$

Define for $(m', n'] = Int_{k-1}((m,n], p_1, \ldots, p_{k-1})$

$$SubFm_k(x, [i,j], (m,n], p_1, \ldots, p_k) =_{df} [c_k + 1, d_k]$$
$$= \bigcap_{1 \leq l \leq k} SubFm_1(x, [i,j], Int_{l-1}((m,n], p_1, \ldots, p_{l-1}), p_l)$$
$$= SubFm_{k-1}(x, [i,j], (m,n], p_1, \ldots, p_{k-1}) \cap SubFm_1(x, [i,j], (m', n'], p_k)$$

**Lemma 4.2** (=Lemma 13 in p.22 [7].) *AID proves the following: Suppose $x[i,j]$ is a $\leq 1$-scarred subformula, $n < i \leq j \leq m$, $n - m = \Delta_{u+1}$, $k \geq 0$ and $1 \leq p_1, \ldots, p_k \leq 4$, and let $A$ denote the interval $SubFm_k(x, [i,j], (m,n], p_1, \ldots, p_k)$. Then*

(a) $A$ is properly contained in $Int_k((m,n), p_1, \ldots, p_k)$.

(b) Each symbol in $A$ is in exactly one of the intervals
$SubFm_{k+1}(x, [i,j], (m,n), p_1, \ldots, p_k, p_{k+1})$ $(1 \leq p_{k+1} \leq 4)$ or is the binary operator $BinOp(x, A, Int_k((m,n), p_1, \ldots, p_k))$.

(c) Each $SubFm_{k+1}(x, [i,j], (m,n), p_1, \ldots, p_k, p_{k+1})$ $(1 \leq p_{k+1} \leq 4)$ is a $\leq 1$-scarred subformula.

Finally define the truth value of $x[i,j]$ by synthesizing truth values $Value_k(x, [i,j], (m,n), \vec{p})$ of subformulae $SubFm_k(x, [i,j], (m,n), \vec{p})$ $(\vec{p} = p_1, \ldots, p_k)$.

**Definition 4.2** in p.23 [7]. Let $n - m = \Delta_{u+1}, m < i \leq j \leq n, x[i,j]$ is a $\leq 1$-scarred formula, $0 \leq k \leq u + 1, 1 \leq p_1, \ldots, p_k \leq 4$.

$Value_k(x, [i,j], (m,n), \vec{p})$ $(\vec{p} = p_1, \ldots, p_k)$ is defined by:

**Case 1** $k = u + 1$: Then $Int_k(x, (m,n), \vec{p}) = (a, a+2]$ for some $a$.
If $SubFm_k(x, [i,j], (m,n), \vec{p})$ is undefined, then

$$Value_k(x, [i,j], (m,n), \vec{p}) =_{df} (\top, \bot)$$

Otherwise $SubFm_k(x, [i,j], (m,n), \vec{p})$ consists of a single logival symbol, $\neg$ or $\top$ or $\bot$ or a variable $q$. Then
$Value_k(x, [i,j], (m,n), \vec{p})$ is defined to be $(\bot, \top)$ or $(\top, \top)$ or $(\bot, \bot)$ or $(q, q)$, resp. By $(q, q)$ we mean $(\top, \top)$ if $q$ has truth value *True* and $(\bot, \bot)$ if $q$ has truth value *False*.

**Cases 2 and 3** $k \leq u$: Let $(a, b] = Int_k(x, [i,j], (m,n), \vec{p})$, $1 \leq p_{k+1} \leq 4$. Put $I_{p_{k+1}} = SubFm_{k+1}(x, [i,j], (m,n), \vec{p}, p_{k+1})$. Then by Lemma 4.2.(b)
$A = SubFm_k(x, [i,j], (m,n), \vec{p}) = I_1 \dot\cup I_2 \dot\cup I_3 \dot\cup I_4 \dot\cup BinOp(x, A, Int_k((m,n), \vec{p}))$
(disjoin union) by using breakpoints $a_1, \ldots, a_4$ of $x[i,j]$ generated by $(a, b]$. Let $v_p$ $(p = 1, 2, 3, 4)$ denote the truth value

$$v_p = Value_{k+1}(x, [i,j], (m,n), \vec{p}, p)$$

**Case 2** $a_2 = a_4$: Then define

$$Value_k(x, [i,j], (m,n), \vec{p}) = v_1 \circ v_2 \circ v_4$$

for the (reverse) composition $\circ$:

$$(r_1, r_2) = (s_1, s_2) \circ (t_1, t_2) \Leftrightarrow_{df} r_i = \begin{cases} t_1 & \text{if } s_i = \top \\ t_2 & \text{if } s_i = \bot \end{cases}$$

**Case 3** $a_2 \neq a_4$: Then define

$$Value_k(x, [i,j], (m,n), \vec{p}) = f_{BinOp}(v_1, v_2) \circ v_3 \circ v_4$$

where, if $BinOp(x, [i,j], Int_k((m,n), \vec{p})) = \odot$, then

$$f_{BinOp}((s_1, s_2), (t_1, t_2)) =_{df} (s_1 \odot t_1, s_2 \odot t_2)$$

(In this case we have $t_1 = t_2$.)

Now we examine the definition of $Value_k(x,[i,j],(m,n),p_1,\ldots,p_k)$ in $AID$.
Let $Value_k^\xi(x,[i,j],(m,n),\vec{p})$ ($\vec{p} = p_1,\ldots,p_k$) be a predicate for $\xi \in \{\top,\bot\}$ such that if
$$Value_k(x,[i,j],(m,n),\vec{p}) = (\xi_1,\xi_2)\ (\xi_1,\xi_2 \in \{\top,\bot\}),$$
then
$$Value_k^\top(x,[i,j],(m,n),\vec{p})\ \text{holds} \Leftrightarrow \xi_1 = \top$$
$$Value_k^\bot(x,[i,j],(m,n),\vec{p})\ \text{holds} \Leftrightarrow \xi_2 = \top$$

Definition 4.2 gives a *simultaneous inductive definition* of the predicates $Value_k^\xi(x,[i,j],(m,n),\vec{p})$ ($\xi \in \{\top,\bot\}$) along a quadtree of depth $u+2 < \|j-i\|+2 \leq \|x\|_\Sigma\|+2 \leq \|x\|+2$: for some $\Sigma_0^b$ $B, \bar{D}$ in $L_{AID}$ and a boolean $I$.

**Case 1** $k = u+1$: $Value_k^\xi(x,[i,j],(m,n),\vec{p}) \leftrightarrow B(x,i,j,m,n,\vec{p},\xi)$

**Cases 2 and 3** $k \leq u$: $Value_k^\xi(x,[i,j],(m,n),\vec{p})$ iff
$I(\bar{D}(x,i,j,m,n,\vec{p},\xi),\{Value_{k+1}^\eta(x,[i,j],(m,n),\vec{p},i) : i = 1,\ldots,4, \eta \in \{\top,\bot\}\})$
($BinOp(x,[i,j],Int_k((m,n),\vec{p}))$ is one of the finitely many binary connectives, and so can be written as a finite disjunction.)

Thus by Lemmata 2.3, 2.4 and 2.5, $Value_k^\xi(x,[i,j],(m,n),\vec{p})$ is $\Sigma_0^b$-definable in $AID$.

Further this truth definition respects the meanings of propositional connectives and the truth value $Value_0(x,[i,j],(m,n))$ is independent of $m$ and $n$.

Assume that $x[i,j]$ is a $\leq 1$-scarred formula, $n - m = \Delta_{u+1}, m < i \leq j \leq n$. Then we have the following lemmata and corollary by $\Sigma_0^b - LIND$.

**Lemma 4.3** (=Lemma 14 in p.24 [7].) *If $x[j]$ is a unary connective $\odot$, then $Value_0(x,[i,j],(m,n)) = Value_0(x,[i,j-1],(m,n)) \circ (s_1,s_2)$ where $(s_1,s_2)$ is the pair of boolean truth values giving the truth value of the $\leq 1$-scarred formula $\odot$. (If $\odot = \neg$, then $(s_1,s_2) = (\bot,\top)$)*

**Lemma 4.4** (=Lemma 15 in p.25 [7].) *Suppose that $x[j]$ is a binary connective $\odot$ and let $f_\odot$ be the binary function such that $f_\odot((s_1,s_2),(t,t)) = (s_1 \odot t, s_2 \odot t)$. Then*

**(a)** *If $x[i,j-1]$ is an unscarred formula,*

$$Value_0(x,[i,j],(m,n)) = f_\odot((\top,\bot), Value_0(x,[i,j-1],(m,n)))$$

**(b)** *Otherwise, let $k \in [i,j]$ be such that $x[k,j-1]$ is a formula (unscarred). Then*

$$Value_0(x,[i,j],(m,n)) =$$
$$f_\odot(Value_0(x,[i,k-1],(m,n)), Value_0(x,[k,j-1],(m,n)))$$

**Corollary 4.1** (=Corollary 16 in p.25 [7].) *If $x[i, j]$ is a formula and if $m_k < i \leq j \leq j_k, n_k - m_k = \Delta_{u_k+1}$ for $k = 1, 2$, then*

$$Value_0(x, [i, j], (m_1, n_1)) = Value_0(x, [i, j], (m_2, n_2))$$

Thus we can $\Sigma_0^b$-define the truth for PLOF formulae by

$$TRUE_{PLOF}(x, [i, j]) = Value_0^\xi(x, [i, j], (m, n))$$

with $\xi \in \{\top, \bot\}, n-m = \Delta_{u+1}$ and $m < i \leq j \leq n$, e.g., $m = i-1, u = \|j-i\|-1$.
Let $RFN(PLOF - \mathcal{F})$ denote

$$\forall x \in PLOF[PLOF - \mathcal{F} \vdash x \to TRUE_{PLOF}(x)],$$

Reflection schema for a PLOF Frege system. Proofs in $PLOF - \mathcal{F}$ are sequences of PLOF formulae separated commas. By counting commas we can $\Sigma_0^b$-define a function $\beta(i, x) = i^{th}$ formula of a proof $x$.

For an infix formula $x$ we set

$$TRUE(x) \Leftrightarrow_{df} TRUE_{PLOF}(\{i < |x| : PL(x, i)\})$$

for the PLOF form $y = \{i < |x| : PL(x, i)\}$ of $x$, cf. Proposition 4.1, Definitions 5.5, 5.6 and Lemma 5.10. Thus $RFN(\mathcal{F})$ denotes $\forall x[\mathcal{F} \vdash x \to TRUE(x)]$, i.e., $\forall x[\mathcal{F} \vdash x \to TRUE_{PLOF}(\{i < |x| : PL(x, i)\})]$, Reflection schema for a Frege system.

**Theorem 4.1**     1. $AID \vdash RFN(PLOF - \mathcal{F})$ *for any Frege system* $\mathcal{F}$.

  2. $AID \vdash RFN(\mathcal{F})$ *for any Frege system* $\mathcal{F}$.

# 5 Stratifications

In this section *stratifications* of formulae are defined. These are in essence to interprete first order formulae into second order formulae. In later sections we need these.

First we define stratified formulae in $L_{BA}$.

**Definition 5.1** (Stratified formulae in $L_{BA}$)
Let $\bar{x}$ and $\bar{i}$ be sequences of variables with $\bar{x} \cap \bar{i} = \emptyset$. Let $\mathcal{B}_{BA}(\bar{x}; \bar{i})$ denote the set of formulae generated as follows:

  1. Atomic formulae of the forms

$$Bit(\ell(\bar{i}), x) = 1, Bit(\ell(\bar{i}), |x|) = 1, Bit(\ell(\bar{i}), |x| \cdot |y|) = 1,$$

$$Bit(\ell(\bar{i}), i) = 1, Bit(\ell(\bar{i}), |i|) = 1, i < |x|, \ell(\bar{i}) \leq \ell'(\bar{i})$$

are in $\mathcal{B}_{BA}(\bar{x}; \bar{i})$, where $x, y$ are in the list $\bar{x}$, $i$ is in the list $\bar{i}$, $\ell(\bar{i}), \ell'(\bar{i})$ are linear forms $\ell(\bar{i}) = \sum_k c_k i_k + d$ with $\bar{i} = i_1, \ldots$

2. If $B_0(\bar{x};\bar{\imath}), B_1(\bar{x};\bar{\imath}) \in \mathcal{B}_{BA}(\bar{x};\bar{\imath})$, then $B_0 \wedge B_1, B_0 \vee B_1, \neg B_0 \in \mathcal{B}_{BA}(\bar{x};\bar{\imath})$.

3. If $B(\bar{x};\bar{\imath}^\frown j) \in \mathcal{B}_{BA}(\bar{x};\bar{\imath}^\frown j)$, then $Qj < p|\bar{x}| \, B(\bar{x};\bar{\imath}^\frown j) \in \mathcal{B}_{BA}(\bar{x};\bar{\imath})$ for any polynomial $p|\bar{x}|$ and $Q \in \{\forall, \exists\}$.

4. If $B(\bar{x}^\frown y;\bar{\imath}) \in \mathcal{B}_{BA}(\bar{x}^\frown y;\bar{\imath})$, then $Q|y| \leq \ell\|\bar{x}\| B(\bar{x}^\frown y;\bar{\imath}) \in \mathcal{B}_{BA}(\bar{x};\bar{\imath})$ for any linear form $\ell\|\bar{x}\|$ and $Q \in \{\forall, \exists\}$.

We say that a $\Sigma_0^b$-formula $B(\bar{x};\bar{\imath})$ in $L_{BA}$ is *stratified with respect to* $(\bar{x};\bar{\imath})$ if $B(\bar{x};\bar{\imath}) \in \mathcal{B}_{BA}(\bar{x};\bar{\imath})$. Also we say that a $\Sigma_0^b$-formula $B(\bar{x})$ in $L_{BA}$ is *stratified with respect to* $\bar{x}$ if $B(\bar{x}) \in \mathcal{B}_{BA}(\bar{x};)$ with the empty list $\bar{\imath} = \emptyset$. When no confusion likely occurs, we simply write $\mathcal{B}_{BA}$ for $\mathcal{B}_{BA}(\bar{x};\bar{\imath})$.

Observe that function 'constants' occurring in a formula in $\mathcal{B}_{BA}$ are $Bit, +, 0, 1, |\cdot|$ and $|x| \cdot |y|$.

**Definition 5.2** (Bitwise computability)

1. Let $t(\bar{x})$ be a term with variables $\bar{x}$. (Every variable in $t(\bar{x})$ need not be in the list $\bar{x}$.) We say that $t(\bar{x})$ is *bitwise computable with respect to* $\bar{x}$ denoted by $t(\bar{x}) \in \mathcal{BC}_{BA}$ if there exists a stratified formula $C_t^*(\bar{x};i) \in \mathcal{B}_{BA}(\bar{x};i)$ such that
$$AID \vdash \forall i < |t(\bar{x})| \, [Bit(i, t(\bar{x})) = 1 \leftrightarrow C_t^*(\bar{x};i)] \tag{7}$$

2. A term $t(\bar{x})$ is said to be *hereditarily bitwise computable with respect to* $\bar{x}$ denoted by $t(\bar{x}) \in \bar{\mathcal{C}}_{BA}$ if $s(\bar{x}) \in \mathcal{B}_{BA}$ for every subterm $s(\bar{x})$ of $t(\bar{x})$.

3. Let $C(\bar{x})$ be a $\Sigma_0^b$-formula in $L_{BA}$ with variables $\bar{x}$. (Every free variable in $C(\bar{x})$ need not be in the list $\bar{x}$.) We say that $C(\bar{x})$ is *bitwise computable with respect to* $\bar{x}$ denoted by $C(\bar{x}) \in \mathcal{C}_{BA}$ if there exists a stratified formula $C^*(\bar{x}) \in \mathcal{B}_{BA}(\bar{x};)$ such that
$$AID \vdash C(\bar{x}) \leftrightarrow C^*(\bar{x})$$

The following Lemmata 5.1-5.5 are preparatory steps towards showing Lemma 5.6: for any term $t(\bar{x})$, $t(\bar{x}) \in \mathcal{C}_{BA}$.

**Lemma 5.1** $t(\bar{x}), s(\bar{x}) \in \mathcal{C}_{BA} \Rightarrow t(\bar{x}) = s(\bar{x}), t(\bar{x}) < s(\bar{x}) \in \mathcal{C}_{BA}$

**Proof**. This follows from the facts: for any $z$ with $\max\{|x|, |y|\} \leq z$
$$x = y \leftrightarrow \forall i < z(Bit(i,x) = 1 \leftrightarrow Bit(i,y) = 1),$$

$x < y \leftrightarrow$
$$\exists j < z[Bit(j,x) \neq 1 \wedge Bit(j,y) = 1 \wedge$$
$$\forall i < z(j < i \rightarrow (Bit(i,x) = 1 \leftrightarrow Bit(i,y) = 1))],$$

and $\max\{|t(\bar{x})|, |s(\bar{x})|\} \leq p|\bar{x}|$ for a polynomial $p|\bar{x}|$. □

**Lemma 5.2** *For each function constant $f \notin \{|x|, x \cdot 2^{|y|}, x\#y\}$ in $L_{BA}$*

$$\bar{t}(\bar{x}) \in \mathcal{C}_{BA} \Rightarrow f(\bar{t}(\bar{x})) \in \mathcal{C}_{BA}$$

**Proof.**

1. (zero), (one): $Bit(i, 0) = 1 \leftrightarrow \bot$ and $Bit(i, 1) = 1 \leftrightarrow i = 0$.

2. (modified subtraction): $Bit(i, x \dot{-} y) = 1$ iff
    $y < x \wedge i < |x| \wedge [\{(Bit(i,x) = 1 \oplus Bit(i,y) = 1) \wedge y[0,i) \leq x[0,i)\} \vee \{Bit(i,x) = Bit(i,y) \wedge x[0,i) < y[0,i)\}]$ and
    $s(\bar{x})[0,i) \leq t(\bar{x})[0,i), t(\bar{x})[0,i) < s(\bar{x})[0,i) \in \mathcal{C}_{BA}(\bar{x}; i)$ if $t(\bar{x}), s(\bar{x}) \in \mathcal{C}_{BA}$.

3. (half): $Bit(i, \lfloor \frac{x}{2} \rfloor) = 1 \leftrightarrow Bit(i+1, x) = 1$

4. (addition): $Bit(i, x + y) = 1$ iff
    $Bit(i, x) = 1 \oplus Bit(i, y) = 1 \oplus$
    $\exists j < |x| [j < i \wedge Bit(j, x) = 1 \wedge Bit(j, y) = 1 \wedge \forall k < |x| \{j < k < i \rightarrow (Bit(k, x) = 1 \oplus Bit(k, y) = 1)\}] \wedge (i \leq |x| \vee i \leq |y|)$.

5. (part): $Bit(i, x[y, z)) = 1 \leftrightarrow \exists u < |x| (Bit(u, x) = 1 \wedge y + i = u < z)$

□

**Lemma 5.3** $t(\bar{x}) \in \bar{\mathcal{C}}_{BA} \Rightarrow |t(\bar{x})| \in \mathcal{C}_{BA}$

**Proof** by induction on the length of the term $t(\bar{x}) \in \bar{\mathcal{C}}_{BA}$.

1. (zero), (one): $Bit(i, |0|) = 1 \leftrightarrow \bot$ and $Bit(i, |1|) = 1 \leftrightarrow i = 0$.

2. (variable): $C^*_{|x|}(x; i) \leftrightarrow_{df} Bit(i, |x|) = 1$.

3. (length): $Bit(i, \|t(\bar{x})\|) = 1 \leftrightarrow \exists |y| \leq \ell\|\bar{x}\|(y = |t(\bar{x})| \wedge Bit(i, |y|) = 1)$ for a linear form $\ell$ such that $\|t(\bar{x})\| \leq \ell\|\bar{x}\|$.

4. (modified subtraction):
    $|x \dot{-} y| = 1 + \max\{j < |x| : Bit(j, x \dot{-} y) = 1\}$. Hence $Bit(i, \max\{j < |x| : Bit(j, x \dot{-} y) = 1\}) = 1$ iff $\exists j < |x| [Bit(j, x \dot{-} y) = 1 \wedge Bit(i, j) = 1 \wedge \forall k < j(Bit(k, x \dot{-} y) \neq 1)]$. Then use Lemma 5.2.

5. (half): $Bit(i, |\lfloor \frac{x}{2} \rfloor|) = 1 \leftrightarrow Bit(i, |x| \dot{-} 1) = 1$. The case (modified subtraction) in Lemma 5.2 yields the assertion.

6. (part): By $|x[i,j]| = \min\{|x|, j\} \dot{-} i$ and (modified subtraction) it suffices to show:

    $$t(\bar{x}), s(\bar{x}) \in \bar{\mathcal{C}}_{BA} \Rightarrow \min\{t(\bar{x}), s(\bar{x})\}, \max\{t(\bar{x}), s(\bar{x})\} \in \mathcal{C}_{BA}$$

    These follow from Lemma 5.1 and $Bit(i, \min\{t(\bar{x}), s(\bar{x})\}) = 1$ iff $(t(\bar{x}) \leq s(\bar{x}) \wedge Bit(i, t(\bar{x})) = 1) \vee (s(\bar{x}) < t(\bar{x}) \wedge Bit(i, s(\bar{x})) = 1)$, and similarly for max.

7. (addition): $\max\{|x|,|y|\} \leq |x+y| \leq \max\{|x|,|y|\}+1$ and
   $|x+y| = \max\{|x|,|y|\}+1$ iff $Bit(\max\{|x|,|y|\}, x+y) = 1$ iff
   $\exists u \leq |x|+|y|\ (u = \max\{|x|,|y|\} \wedge Bit(u, x+y) = 1)$. Further
   $(u = \max\{|t(\bar{x})|, |s(\bar{x})|\} \wedge Bit(u, t(\bar{x})+s(\bar{x})) = 1) \in \mathcal{C}_{BA}(\bar{x}; u)$ if $|t(\bar{x})|, |s(\bar{x})| \in \bar{\mathcal{C}}_{BA}$.

8. (smash): $|x \# y| = |x| \cdot |y| + 1$.

9. (padding): $|x \cdot 2^{|y|}| = |x| + |y|$. These are handled by Lemma 5.2.

$\square$

**Lemma 5.4** $t(\bar{x}), s(\bar{x}) \in \bar{\mathcal{C}}_{BA} \Rightarrow t(\bar{x}) \cdot 2^{|s(\bar{x})|} \in \mathcal{C}_{BA}$

**Proof.** This follows from $Bit(i, x \cdot 2^{|y|}) = 1$ iff $\exists j < |x|\ [Bit(j, x) = 1 \wedge i = j + |y|]$ and Lemmata 5.2, 5.3. $\square$

**Lemma 5.5** $t(\bar{x}), s(\bar{x}) \in \bar{\mathcal{C}}_{BA} \Rightarrow t(\bar{x}) \# s(\bar{x}) \in \mathcal{C}_{BA}$

**Proof.** We have $Bit(i, x \# y) = 1 \leftrightarrow i = |x| \cdot |y|$. Thus it suffices to show

$$t(\bar{x}), s(\bar{x}) \in \bar{\mathcal{C}}_{BA} \Rightarrow |t(\bar{x})| \cdot |s(\bar{x})| \in \bar{\mathcal{C}}_{BA}$$

This is seen as in Lemma 5.3 by induction on the sum of the lengths of terms $t(\bar{x}), s(\bar{x})$. Note that by the definition $Bit(i, |x| \cdot |y|) = 1 \in \mathcal{B}_{BA}(x, y; i)$. $\square$

By Lemmata 5.2, 5.3, 5.4 and 5.5 we get

**Lemma 5.6** *For any term $t(\bar{x})$, $t(\bar{x}) \in \mathcal{C}_{BA}$.*

**Definition 5.3** For a term $t(\bar{x})$ $C_t^*(\bar{x}; i)$ denotes a stratified ($\Sigma_0^b$-)formula in $\mathcal{B}_{BA}(\bar{x}; i)$ so that

$$AID \vdash \forall i < |t(\bar{x})|\ [Bit(i, t(\bar{x})) = 1 \leftrightarrow C_t^*(\bar{x}; i)] \tag{7}$$

Next we define stratified formulae in $L_{AID}$.

**Definition 5.4** (Stratified formulae in $L_{AID}$)
Let $\bar{x}$ and $\bar{\imath}$ be sequences of variables with $\bar{x} \cap \bar{\imath} = \emptyset$. Let $\mathcal{B}(\bar{x}; \bar{\imath})$ denote the set of formulae generated as follows:

1. $\mathcal{B}_{BA}(\bar{x}; \bar{\imath}) \subseteq \mathcal{B}(\bar{x}; \bar{\imath})$.

2. For the inductively defined predicate $A^{\ell, B, \bar{D}, I}$ defined from $\ell, B, \bar{D}, I$, $A^{\ell, B, \bar{D}, I}(t_1, \ldots, t_n, s) \in \mathcal{B}(\bar{x}; \bar{\imath})$ iff

   (a) Terms $t_1, \ldots, t_n, s$ are variables $y_1, \ldots, y_n, p$ so that $y_1, \ldots, y_n \subseteq \bar{x}$ & $p \in \bar{\imath}$, and

   (b) $B(x_1, \ldots, x_n, i), \bar{D}(x_1, \ldots, x_n, i) \in \mathcal{B}_{BA}(\bar{x}; \bar{\imath})$.

3. $\mathcal{B}(\bar{x};\bar{\imath})$ is closed under propositional connectives.

4. If $B(\bar{x};\bar{\imath}\frown j) \in \mathcal{B}(\bar{x};\bar{\imath}\frown j)$ , then $Qj < p\,|\,\bar{x}\,|\, B(\bar{x};\bar{\imath}\frown j) \in \mathcal{B}(\bar{x};\bar{\imath})$ for any polynomial $p\,|\,\bar{x}\,|$ and $Q \in \{\forall, \exists\}$.

5. If $B(\bar{x}\frown y;\bar{\imath}) \in \mathcal{B}(\bar{x}\frown y;\bar{\imath})$ , then $Q\,|\,y\,| \leq \ell \|\bar{x}\| B(\bar{x}\frown y;\bar{\imath}) \in \mathcal{B}(\bar{x};\bar{\imath})$ for any linear form $\ell \|\bar{x}\|$ and $Q \in \{\forall, \exists\}$.

We say that a $\Sigma_0^b$-formula $B(\bar{x};\bar{\imath})$ in $L_{AID}$ is *stratified with respect to* $(\bar{x};\bar{\imath})$ if $B(\bar{x};\bar{\imath}) \in \mathcal{B}(\bar{x};\bar{\imath})$. Also we say that a $\Sigma_0^b$-formula $B(\bar{x})$ in $L_{AID}$ is *stratified with respect to* $\bar{x}$ if $B(\bar{x}) \in \mathcal{B}(\bar{x};)$ with the empty list $\bar{\imath} = \emptyset$. When no confusion likely occurs, we simply write $\mathcal{B}$ for $\mathcal{B}(\bar{x};\bar{\imath})$.

Let $C(\bar{x})$ be a $\Sigma_0^b$-formula in $L_{AID}$ with variables $\bar{x}$. (Every free variable in $C(\bar{x})$ need not be in the list $\bar{x}$.) We say that $C(\bar{x})$ is *bitwise computable with respect to* $\bar{x}$ denoted by $C(\bar{x}) \in \mathcal{C}$ if there exists a stratified formula $C^*(\bar{x}) \in \mathcal{B}(\bar{x};)$ such that
$$AID \vdash C(\bar{x}) \leftrightarrow C^*(\bar{x})$$

In what follows we show that Lemma 5.9, $C(\bar{x}) \in \mathcal{C}$ for any $\Sigma_0^b$ $C$.

**Definition 5.5** (Substituting a formula in a stratified formula)
Let $C(y) \in \mathcal{B}(\bar{y};\bar{\jmath})$ with $y \in \bar{y}$ and $A_0(\bar{x}, i)$ a $\Sigma_0^b$-formula and $p\,|\,\bar{x}\,|$ a polynomial. We define a $\Sigma_0^b$-formula $C(\{i < p\,|\,\bar{x}\,|: A_0(\bar{x}, i)\})$ as follows: Let $lh(p, A_0)$ denote a $\Sigma_0^b$-definable function in $AID$ such that

$$lh(p, A_0) = \begin{cases} \max\{i < p\,|\,\bar{x}\,|: A_0(\bar{x}, i)\} + 1 & \text{if } \exists i < p\,|\,\bar{x}\,|\, A_0(\bar{x}, i) \\ 0 & \text{otherwise} \end{cases}$$

Then $C(\{i < p\,|\,\bar{x}\,|: A_0(\bar{x}, i)\})$ is obtained by replacing $Bit(j, y) = 1$ by $j < p\,|\,\bar{x}\,| \wedge A_0(\bar{x}, j)$ and by replacing $|y|$ by $lh(p, A_0)$.

1. The substitution $y \to \{i < p\,|\,\bar{x}\,|: A_0(\bar{x}, i)\}$ commutes with propositional connectives.

2. If $C(y) \equiv Qz < q(|y|)C_0(z, y)$ , then $C(\{i < p\,|\,\bar{x}\,|: A_0(\bar{x}, i)\}) \Leftrightarrow_{df}$ $Qz < q(lh(p, A))C_0(z, \{i < p\,|\,\bar{x}\,|: A_0(\bar{x}, i)\})$.

3. (a) $Bit(j, \{i < p\,|\,\bar{x}\,|: A_0(\bar{x}, i)\}) = 1 \Leftrightarrow_{df} j < p\,|\,\bar{x}\,| \wedge A_0(\bar{x}, j)$,
    (b) $Bit(j, |\{i < p\,|\,\bar{x}\,|: A_0(\bar{x}, i)\}|) = 1 \Leftrightarrow_{df} Bit(j, lh(p, A_0)) = 1$,
    (c) $Bit(j, |\{i < p\,|\,\bar{x}\,|: A_0(\bar{x}, i)\}| \cdot |y_1|) = 1 \Leftrightarrow_{df}$ $Bit(j, lh(p, A_0) \cdot |y_1|) = 1$, and
    (d) $Bit(j, |\{i < p\,|\,\bar{x}\,|: A_0(\bar{x}, i)\}| \cdot |\{i < p\,|\,\bar{x}\,|: A_0(\bar{x}, i)\}|) = 1 \Leftrightarrow_{df}$ $Bit(j, lh(p, A_0) \cdot lh(p, A_0)) = 1$.

4. $j < |\{i < p\,|\,\bar{x}\,|: A_0(\bar{x}, i)\}| \Leftrightarrow_{df} j < lh(p, A_0)$

5. The case $C(y) \equiv A^{\ell, B, \bar{D}, I}(t_1, \ldots, t_{i-1}, y, t_{i+1}, \ldots, t_n, j)$: For simplicity we assume that none of variables $t_1, \ldots, t_{i-1}, t_{i+1}, \ldots, t_n$ is the variable $y$. Let $\ell'$ denote a linear form such that if $\ell \|\bar{z}\| = \sum_k c_k \|z_k\| + d$, then

$$\sum_{k \neq i} c_k \|z_k\| + d + c_i \,|lh(p, A_0)| < \ell' \|z_1, \ldots, z_{i-1}, z_{i+1}, \ldots, z_n, \bar{x}\|$$

Put

$$B'(z_1, \ldots, z_{i-1}, z_{i+1}, \ldots, z_n, \bar{x}, p) \Leftrightarrow_{df}$$
$$B(z_1, \ldots, z_{i-1}, \{i < p \,|\bar{x}|: A_0(\bar{x}, i)\}, z_{i+1}, \ldots, z_n, \bar{x}, p)$$
$$\bar{D}'(z_1, \ldots, z_{i-1}, z_{i+1}, \ldots, z_n, \bar{x}, p) \Leftrightarrow_{df}$$
$$\bar{D}(z_1, \ldots, z_{i-1}, \{i < p \,|\bar{x}|: A_0(\bar{x}, i)\}, z_{i+1}, \ldots, z_n, \bar{x}, p)$$

$B', \bar{D}'$ are $\Sigma_0^b$-formuale in $L_{AID}$ if $A_0(\bar{x}, i) \notin L_{BA}$.
Let $A'(z_1, \ldots, z_{i-1}, z_{i+1}, \ldots, z_n, \bar{x}, p)$ denote the inductively defined predicate such that for $0 \neq |p| \leq \ell' \|z_1, \ldots, z_{i-1}, z_{i+1}, \ldots, z_n, \bar{x}\|$
$A'(z_1, \ldots, z_{i-1}, z_{i+1}, \ldots, z_n, \bar{x}, p)$ iff
$|p| = \sum_{k \neq i} c_k \|z_k\| + d + c_i \,|lh(p, A_0)|$ & $B'(\bar{z}, \bar{x}, p)$ or
$|p| < \sum_{k \neq i} c_k \|z_k\| + d + c_i \,|lh(p, A_0)|$ and
$I(\bar{D}'(\bar{z}, \bar{x}, p), A'(\bar{z}, \bar{x}, p0), A'(\bar{z}, \bar{x}, p1))$
with $\bar{z} = z_1, \ldots, z_{i-1}, z_{i+1}, \ldots, z_n$.
$A'$ is $\Sigma_0^b$-definable in $AID$ by Lemma 2.4, (Iterated inductive definitions) when $A_0 \notin L_{BA}$. Then

$$A^{\ell, B, \bar{D}, I}(t_1, \ldots, t_{i-1}, \{i < p \,|\bar{x}|: A_0(\bar{x}, i)\}, t_{i+1}, \ldots, t_n, j) \Leftrightarrow_{df}$$
$$A'(t_1, \ldots, t_{i-1}, t_{i+1}, \ldots, t_n, \bar{x}, j)$$

**Lemma 5.7** *For a stratified formula $C^*(y)$*

$$AID \vdash y = \{i < p \,|\bar{x}|: A(\bar{x}, i)\} \rightarrow [C^*(y) \leftrightarrow C^*(\{i < p \,|\bar{x}|: A(\bar{x}, i)\})]$$

*where $y = \{i < p \,|\bar{x}|: A(\bar{x}, i)\} \Leftrightarrow_{df} |y| < p \,|\bar{x}|$ & $\forall i < p \,|\bar{x}| \,(i \in y \leftrightarrow A(\bar{x}, i))$.*

**Lemma 5.8** *For any terms $\bar{t}(\bar{x})$, $C(\bar{x}, \bar{y}) \in \mathcal{C} \Rightarrow C(\bar{x}, \bar{t}(\bar{x})) \in \mathcal{C}$.*

**Proof**. Let $C^*(\bar{x}, \bar{y})$ be a stratified formula such that
$AID \vdash C(\bar{x}, \bar{y}) \leftrightarrow C^*(\bar{x}, \bar{y})$. By Lemma 5.6 we have $\bar{t}(\bar{x}) \in \mathcal{C}$. First show the case $C \in L_{BA}$ by induction on $C^*$ and then the general case again by induction on $C^*$. For the case $C^*(y)$ is a formula $A(\ldots, y, \ldots, j)$ with an inductively defined predicate $A$ use the first case and Definition 5.5. □

**Lemma 5.9** *For any $\Sigma_0^b$-formula $C(\bar{x}) \in L_{AID}$, $C(\bar{x}) \in \in \mathcal{C}$.*

**Proof**. By induction on the length of formulae using Lemmata 5.6 and 5.8. If $C(\bar{x})$ is a formula $A(\bar{x}, t(\bar{x}))$ for an inductively defined predicate $A$, then put $A(\bar{x}, t(\bar{x})) \leftrightarrow \exists \mid u \mid \leq \ell \|\bar{x}\| (A(\bar{x}, u) \wedge u = t(\bar{x}))$. □

**Definition 5.6**   *1. For a $\Sigma_0^b$-formula $C(\bar{x})$ $C^*(\bar{x})$ denotes a stratified ($\Sigma_0^b$-)formula in $\mathcal{B}(\bar{x};)$ so that*

$$AID \vdash C(\bar{x}) \leftrightarrow C^*(\bar{x})$$

*2. For $\Sigma_0^b$-formulae $C(y), A_0(\bar{x}, i)$ and a polynomial $p \mid \bar{x} \mid C(\{i < p \mid \bar{x} \mid: A_0(\bar{x}, i)\})$ denotes the formula $C^*(\{i < p \mid \bar{x} \mid: A_0(\bar{x}, i)\})$.*

**Lemma 5.10** *Let $C(y)$ be a $\Sigma_0^b$-formula in $L_{AID}$. Let $f(\bar{x})$ be a $\Sigma_0^b$-bitdefinable function in $AID$, and hence $\Sigma_0^b$-definable in $AID + \Sigma_0^b - CA$:*

$$f(\bar{x}) = y \Leftrightarrow_{df} |y| \leq p \mid \bar{x} \mid \,\&\, \forall i < p \mid \bar{x} \mid (Bit(i, y) = 1 \leftrightarrow A(\bar{x}, i)) \ (A \in \Sigma_0^b)$$

*Then $C(f(\bar{x}))$ is $\Sigma_0^b$-definable in $AID$.*

**Lemma 5.11**   *1. For a term $t(\bar{x})$ let $t^*(\bar{x})$ denote $\{i < |t(\bar{x})| : C_t^*(\bar{x}; i)\}$. For terms $s(z, y)$ and $t(y)$ let $u(y) =_{df} s(t(y), y) =_{df} s(z, y)[t(y)/z]$ denote the result of substituting $t(y)$ for $z$ in $s(z, y)$. Then*

$$AID \vdash u^*(y) = s^*(z, y)[t^*(y)/z]$$

*where $=$ means that these are coextensional.*

*2. For a $\Sigma_0^b$-formula $C(z, y)$ and a term $t(y)$*

$$AID \vdash C(t(y), y) \leftrightarrow C^*(z, y)[t^*(y)/z]$$

The following lemma is neeeded in section 8.

**Lemma 5.12** *Let $B(i)$ be a $\Sigma_0^b$-formula in which a variable $y$ does not occur, $C(z, y)$ a $\Sigma_0^b$-formula, $p \mid \bar{t}_0 \mid$ a polynomial for some terms $\bar{t}_0$, $t(y), s(y)$ terms and $\ell$ a linear form. Let $D(y)$ denote the following $\Sigma_0^b$-formula:*

$$D(y) \equiv |t(y)| \leq \ell \|s(y)\| \wedge C(t(y), y) \to \exists \mid z \mid \leq \ell \|s(y)\| C(z, y)$$

*Then $AID \vdash D^*(y)[\{i < p \mid \bar{t}_0 \mid : B(i)\}/y]$.*

**Proof**. By Lemma 5.11 we have $C(t(y), y) \to C^*(z, y)[t^*(y)/z]$. Let $C_1(y)$ be a stratified formula such that $C_1(y) \leftrightarrow C(t(y), y)$. Then
$C_1(B) \to C^*(z, y)[t^*(B)/z, B/y]$ for $B = \{i < p \mid \bar{t}_0 \mid : B(i)\}$ and $t^*(B) = t^*(y)[B/y]$, a $\Sigma_0^b$-formula. By $|t^*(B)| \leq \ell \|s^*(B)\|$ and $\Sigma_0^b - LCA$, Lemma 1.1 we have
$\exists \mid z \mid \leq \ell \|s^*(B)\| \forall i < |t^*(B)| \, [i \in z \leftrightarrow t^*(B)]$. □

# 6   Frege system simulates $AID$

In this section we show that any $\Sigma_0^b$-theorem in $AID$ yields true boolean sentences of which $\mathcal{F}$ has polysize proofs.

For each stratified $\Sigma_0^b$-formula $B(\bar{x};\bar{\imath}) \in \mathcal{B}(\bar{x};\bar{\imath})$ ($\bar{x} = x_1,\ldots,x_n$) in $L_{AID}$ we define a boolean formula $\langle B(\bar{x};\bar{\imath}) \rangle$ and a valuation
$\sigma_{\bar{x}} : \{p_j^k : 1 \leq k \leq n, j < |x_k|\} \to \{\top,\bot\}$ so that

1. Atoms occurring in $\langle B(\bar{x};\bar{\imath}) \rangle$ are among the atoms $\bar{p}^1,\ldots,\bar{p}^n$ where $\bar{p}^k = p_0^k,\ldots,p_{m_k-1}^k$ with $m_k = |x_k|$.

2. The bitgraph of the function $\varphi_B : (\bar{x};\bar{\imath}) \mapsto \langle B(\bar{x};\bar{\imath}) \rangle$ is $\Sigma_0^b$-definable in $AID$, i.e., there exists a $\Sigma_0^b$-formula $A_B(\bar{x};\bar{\imath},j)$ in $L_{AID}$ for each $B$ such that
$$A_B(\bar{x};\bar{\imath},j) \leftrightarrow Bit(j,\langle B(\bar{x};\bar{\imath}) \rangle) = 1$$
where the boolean formula $\langle B(\bar{x};\bar{\imath}) \rangle$ is coded by 0-1 words as in [7].

3. $AID \vdash B(\bar{x};\bar{\imath}) \leftrightarrow TRUE(\langle B(\bar{x};\bar{\imath}) \rangle; \sigma_{\bar{x}})$ where RHS means that the boolean formula $\langle B(\bar{x};\bar{\imath}) \rangle$ is true under the valuation $\sigma_{\bar{x}}$. Alternatively we can define $\langle B(\bar{x};\bar{\imath}) \rangle$ as a sentence which is the result of replacing each atom $p_j^k$ by $\sigma_{\bar{x}}(p_j^k)$.

**Definition 6.1** (Translation into boolean formulae)

1. The valuation $\sigma_{\bar{x}}$ is defined by
$$\sigma_{\bar{x}}(p_j^k) = \begin{cases} \top & \text{if } Bit(j,x_k) = 1 \\ \bot & \text{if } Bit(j,x_k) = 0 \end{cases}$$

2. 
$$\langle Bit(\ell(\bar{\imath}),x_k) = 1 \rangle = \begin{cases} p_j^k & \text{if } \ell(\bar{\imath}) = j < |x_k| \\ \bot & \text{otherwise} \end{cases}$$

$$\langle Bit(\ell(\bar{\imath}),|x_k| \cdot |x_l|) = 1 \rangle = \begin{cases} \top & \text{if } Bit(\ell(\bar{\imath}),|x_k| \cdot |x_l|) = 1 \\ \bot & \text{otherwise} \end{cases}$$

   and similarly for $\langle Bit(\ell(\bar{\imath}),|x_k|) = 1 \rangle$.

   $\langle Bit(\ell(\bar{\imath}),i) = 1 \rangle, \langle Bit(\ell(\bar{\imath}),|i|) = 1 \rangle, \langle i < |x_k| \rangle, \langle \ell(\bar{\imath}) \leq \ell'(\bar{\imath}) \rangle$ are defined to be $\top$ or $\bot$ if the formula is true or false, resp.

3. Inductively defined predicate $A^{\ell,B,\bar{D},I}$.
   **(A.0)** $\neg[0 \neq |p| \leq \ell\|\bar{x}\|]$: $\langle A^{\ell,B,\bar{D},I}(\bar{x},p) \rangle =_{df} \bot$.
   **(A.1)** $0 \neq |p| \leq \ell\|\bar{x}\|$ & $|p| = \ell\|\bar{x}\|$: $\langle A^{\ell,B,\bar{D},I}(\bar{x},p) \rangle =_{df} \langle B(\bar{x},p) \rangle$
   **(A.2)** $0 \neq |p| \leq \ell\|\bar{x}\|$ & $|p| < \ell\|\bar{x}\|$:
   $\langle A^{B,\bar{D},I,\ell}(\bar{x},p) \rangle =_{df} \langle I(\bar{D}(\bar{x},p), A^{\ell,B,\bar{D},I}(\bar{x},p0), A^{\ell,B,\bar{D},I}(\,barx,p1)) \rangle$.

4. $\langle \cdot \rangle$ commutes with propositional connectives.

5. $\langle \forall \mid y \mid \leq \ell \|\bar{x}\| B(\bar{x}^\frown y; \bar{\imath}) \rangle = \bigwedge \{ \langle B(\bar{x}^\frown y; \bar{\imath}) \rangle : \mid y \mid \leq \ell \|\bar{x}\| \}$ and
$\langle \exists \mid y \mid \leq \ell \|\bar{x}\| B(\bar{x}^\frown y; \bar{\imath}) \rangle = \bigvee \{ \langle B(\bar{x}^\frown y; \bar{\imath}) \rangle : \mid y \mid \leq \ell \|\bar{x}\| \}$.

$\langle Qj < p \mid \bar{x} \mid B(\bar{x}; \bar{\imath}^\frown j) \rangle$ is similarly defined for $Q \in \{\forall, \exists\}$.

It is straightforward to see the

**Theorem 6.1** *For any $\Sigma_0^b$-formula $B(\bar{x})$, if $AID \vdash B(\bar{x})$, then there exists a polynomial $p \mid \bar{x} \mid$ and a $\Sigma_0^b$-formula $P(\bar{x}, i)$ such that*

$$AID \vdash \{i < p \mid \bar{x} \mid : P(\bar{x}, i)\} \text{ is a Frege proof of } \langle B^*(\bar{x}) \rangle$$

*and hence by $\Sigma_0^b - CA$, $\exists \mid y \mid \leq p \mid \bar{x} \mid \forall i < p \mid \bar{x} \mid (i \in y \leftrightarrow P(\bar{x}, i))$, $AID + \Sigma_0^b - CA \vdash \mathcal{F} \vdash^{p|\bar{x}|} \langle B^*(\bar{x}) \rangle$ for an equivalent stratified formula $B^*$.*

Note that $\Sigma_0^b - CA$ is needed here only because of a hidden existential bounded quatifier in $\mathcal{F} \vdash^{p|\bar{x}|}$.

# 7 Theories of bounded arithmetic for Frege

In this section we introduce some systems of bounded arithmetic in the language $L_{BA}$, i.e., without inductively defined predicate $A^{\ell, B, \bar{D}, I}$ which are equivalent to $AID$. These systems contain a base fragment $\Sigma_0^b - LIND$ of bounded arithmetic in the language $L_{BA}$.

**Definition 7.1** $\Sigma_0^b - RD$ *($\Sigma_0^b$-Recursive Definitions) denotes the axiom schema whose instances are of the following form : for $\Sigma_0^b$ $B, \bar{D}$, a boolean $I$ and a linear form $\ell$, cf. **(A.1)** and **(A.2)** in Definition 1.1,*

$$\forall \bar{x} \exists \mid y \mid \leq 2^{\ell\|\bar{x}\|} \forall \mid i \mid \leq \ell\|\bar{x}\| [\{0 \neq \mid i \mid = \ell\|\bar{x}\| \to (i \in y \leftrightarrow B(\bar{x}, i))\} \\ \wedge \{0 \neq \mid i \mid < \ell\|\bar{x}\| \to (i \in y \leftrightarrow I(\bar{D}(\bar{x}, i), i0 \in y, i1 \in y))\}] \tag{8}$$

*where $i \in y \leftrightarrow Bit(i, y) = 1$.*

**Lemma 7.1** $\Sigma_0^b - RD \vdash \Sigma_0^b - CA$

**Proof**. Let $B$ be a $\Sigma_0^b$-formula and $p \mid x \mid$ a polynomial. Let $\ell$ denote a linear form such that $\mid p \mid x \mid \mid < \ell \|x\|$. Pick a $y$ by using $\Sigma_0^b - RD$ so that $\mid i \mid = \ell\|x\| \to (i \in y \leftrightarrow B(i_0) \wedge i_0 < p \mid x \mid)$ for $i_0 = i[0, \ell\|x\| - 1)$. Then $z = y[2^{\ell\|x\|-1}, \mid y \mid)$ is a required set since $i \in z \leftrightarrow 2^{\ell\|x\|-1} + i \in y \leftrightarrow B(i)$ for $i < p \mid x \mid$. $\square$

**Lemma 7.2** *For each inductively defined predicate $A = A^{\ell, B, \bar{D}, I}$ in $L_{AID}$, there exists a $\Delta_1^b$-formula $A'$ in $\Sigma_0^b - RD$, i.e., there are a $\Sigma_1^b$ $A_\Sigma$ and a $\Pi_1^b$ $A_\Pi$ such that $\Sigma_0^b - RD \vdash A' \leftrightarrow_{df} A_\Sigma \leftrightarrow A_\Pi$, so that for any formula $\varphi(A, \ldots)$ in $L_{AID}$,*

$$AID \vdash \varphi(A, \ldots) \Rightarrow \Sigma_0^b - RD + \Delta_1^b - LIND \vdash \varphi(A', \ldots),$$

where each $A$ is replaced by the corresponding $\Delta_1^b$ $A'$.

Note that $\Sigma_0^b - RD + \Delta_1^b - LIND \subseteq \Sigma_0^b - RD + \Delta_1^b - CA \subseteq \Sigma_0^b - RD + \Sigma_1^b - AC$. Thus $AID + \Delta_1^b - CA$ and $AID + \Sigma_1^b - AC$ is interpretable in $\Sigma_0^b - RD + \Delta_1^b - CA$ and $\Sigma_0^b - RD + \Sigma_1^b - AC$, resp.

**Proof**. We show $A = A^{\ell, B, \bar{D}, I}$ is $\Delta_1^b$-definable in $\Sigma_0^b - RD$. Then $\Sigma_0^b(L_{AID}) - LIND$ in $L_{AID}$ turns into $\Delta_1^b(L_{BA}) - LIND$.

Let $p\,|\,\bar{x}\,|$ be a polynomial such that $2^{\ell \|\bar{x}\|} \leq p\,|\,\bar{x}\,|$. Let $Demo(y, \bar{x})$ denote the following $\Sigma_0^b$-formula in $L_{BA}$, cf. (8):

$$|y| \leq p\,|\,\bar{x}\,|\ \&\ \forall\,|\,i\,|\leq \ell\|\bar{x}\|[\{0 \neq |\,i\,| = \ell\|\bar{x}\| \to (i \in y \leftrightarrow B(\bar{x}, i))\}$$
$$\wedge \{0 \neq |\,i\,| < \ell\|\bar{x}\| \to (i \in y \leftrightarrow I(\bar{D}(\bar{x}, i), i0 \in y, i1 \in y))\}]$$

By $\Sigma_0^b - RD$ we have $\forall \bar{x} \exists\,|y| \leq p\,|\,\bar{x}\,|\,Dem(y, \bar{x})$. From $\Sigma_0^b - LIND$ we see that such a demonstration tree $y$ is unique:

$$Demo(y, \bar{x})\ \&\ Demo(z, \bar{x})\ \&\ 0 \neq |\,i\,| \leq \ell\|\bar{x}\| \to (i \in y \leftrightarrow i \in z)$$

Therefore in $\Sigma_0^b - RD$

$$A_\Sigma(\bar{x}, i) \leftrightarrow_{df} \exists\,|y| \leq p\,|\,\bar{x}\,|\,[Demo(y, \bar{x})\ \&\ i \in y]\ \&\ 0 \neq |\,i\,| \leq \ell\|\bar{x}\|$$
$$\leftrightarrow$$
$$A_\Pi(\bar{x}, i) \leftrightarrow_{df} \forall\,|y| \leq p\,|\,\bar{x}\,|\,[Demo(y, \bar{x}) \to i \in y]\ \&\ 0 \neq |\,i\,| \leq \ell\|\bar{x}\|$$

$\square$

# 8 Realizations of $\Sigma_1^b$-consequences

In this section we show that $\Sigma_1^b$-consequences in $AID + \Sigma_0^b - CA$ or $\Sigma_0^b - RD + \Sigma_1^b - AC$ can be realized by a $\Sigma_0^b$-set $\{i < p\,|\,\bar{x}\,|: A(\bar{x}, i)\}_{\bar{x}}$.

Recall that $C^*$ denotes a stratified formula which is equivalent to a given $\Sigma_0^b$-formula $C$, cf. Definition 5.6.

**Lemma 8.1** Let $C(z, \bar{x})$ be a $\Sigma_0^b$-formula and $p_0\,|\,\bar{x}\,|$ a polynomial. If $AID + \Sigma_0^b - CA \vdash \exists\,|z| \leq p_0\,|\,\bar{x}\,|\,C(z, \bar{x})$, then there exists a $\Sigma_0^b$-formula $A(\bar{x}, i)$ such that $AID \vdash C^*(\{i < p_0\,|\,\bar{x}\,|: A(\bar{x}, i)\}, \bar{x})$. In particular $AID + \Sigma_0^b - CA$ is $\Sigma_0^b$-conservative over $AID$.

**Proof**. This is an analogue to the fact about the subsystem $ACA_0$ of second order arithmetic vs. $PA$. Hence the idea of a proof is to replace a 'set' variable $y$ by its $\Sigma_0^b$-definition $\{i < p\,|\,\bar{x}\,|: B(i)\}$ when an instance $\exists\,|y| \leq p\,|\,\bar{x}\,|\,\forall i < p\,|\,\bar{x}\,|\,(i \in y \leftrightarrow B(i))$ of $\Sigma_0^b - CA$ occurs.

Formulate $AID + \Sigma_0^b - CA$ in Gentzen's sequent calculus. $\Sigma_0^b - LIND$ and $\Sigma_0^b - CA$ are replaced by the following inference rules for $B \in \Sigma_0^b$ and eigenvariables $y$:

$$\frac{B(y), \Gamma \to \Delta, B(y+1)}{B(0), \Gamma \to \Delta, B(|t|)}$$

$$\frac{|y| \leq p \, |\bar{x}|, \forall i < p \, |\bar{x}| \, (i \in y \leftrightarrow B(i)), \Gamma \to \Delta}{\Gamma \to \Delta}$$

First eliminate cuts partially to get a proof in which every sequent is $s\Sigma_1^b$. We have $C(A_0) \vee C(A_1) \to C(A)$, where $A_j = \{i < p \, |\bar{x}|: A_j(\bar{x}, i)\}$ and $A(\bar{x}, i) \Leftrightarrow_{df} (C(A_0) \wedge A_0(\bar{x}, i)) \vee (\neg C(A_0) \wedge A_1(\bar{x}, i))$.

This cares $Contraction : right$. $\exists \leq: right$ is seen from Lemma 5.6.

Assume

$$|y| \leq p_0 \, |t_0(\bar{x})|, \forall i < p_0 \, |t_0| \, (i \in y \leftrightarrow B(i)), \Gamma \to \Delta, C(V, \bar{x}) \qquad (9)$$

with $V = \{i < p \, | \bar{x}, y \, |: A(\bar{x}, y, i)\}$. Replace each formula $C_0(y)$ in a proof of the sequent (9) by $C_0^*(\{i < p_0 \, |t_0(\bar{x})| : B(i)\})$. Use Lemma 5.12 to handle $\exists \leq: right$, $\forall \leq: left$ and $\Sigma_0^b - LIND$. Lemma 5.11.2 cares BASIC. Therefore $AID$ proves $\Gamma \to \Delta, C(\{i < p' \, |\bar{x}|: A'(\bar{x}, i)\}, \bar{x})$ with $A'(\bar{x}, i) \Leftrightarrow_{df} A^*(\bar{x}, \{i < p_0 \, |t_0(\bar{x})|: B(i)\}, i)$ and a polynomial $p'$. □

**Lemma 8.2** *For a $\Sigma_0^b$-formula $B(y, \bar{x})$ and a polynomial $p \, |\bar{x}|$, if $\Sigma_0^b - RD + \Sigma_1^b - AC \vdash \forall \bar{x} \exists \, |y| \leq p \, |\bar{x}| \, B(y, \bar{x})$, then there exists a $\Sigma_0^b$-formula $A(\bar{x}, i)$ in $L_{AID}$ such that $AID \vdash B^*(\{i < p \, |\bar{x}|: A(\bar{x}, i)\}, \bar{x})$*

**Proof**. Formulate the system $\Sigma_0^b - RD + \Sigma_1^b - AC = \Sigma_0^b - RD + \Sigma_0^b - AC$ in Gentzen's sequent calculus:

1. Initial sequents: logical ones $A \to A$ $(A \in \Sigma_0^b)$, $BASIC$, Bit Extensionality Axiom, $\Sigma_0^b - LIND$, $\Sigma_0^b - RD$. These are in $s\Sigma_1^b$.

2. Inference rules $LKB$ and $\Sigma_0^b - AC$: for $B \in \Sigma_0^b$

$$\frac{i < p \, |t|, \Gamma \to \Delta, \exists \, |y| \leq q \, |t| \, B(i, y)}{\Gamma \to \Delta, \exists \, |z| \leq p \, |t| \cdot q \, |t| \, \forall i < p \, |t| \, B(i, z_i)}$$

Suppose $\Sigma_0^b - RD + \Sigma_1^b - AC \vdash \forall \bar{x} \exists \, |y| \leq p \, |\bar{x}| \, B(y, \bar{x})$. Eliminate cuts partially. There is a proof of $\to \exists \, |y| \leq p \, |\bar{x}| \, B(y, \bar{x})$ such that every formula in it is either $\Sigma_0^b$ or $s\Sigma_1^b$. That is to say, every sequent in it is of the form:

$$\{\exists \, |z_i| \leq q_i \, |\bar{t}_i| \, C_i(z_i, \bar{x}) : i < n\}, \Pi \to \Lambda, \{\exists \, |y_j| \leq p_j \, |\bar{s}_j| \, B_j(y_j, \bar{x}) : j < m\}$$

with $\{C_i\}_{i<n}, \Pi, \Lambda, \{B_j\}_{j<m} \subseteq \Sigma_0^b$.

We show, by induction on the depth of proofs, that there exist $\Sigma^b_0$-formulae $A_j(\bar{x}, \bar{z}, i)$ $(j < m, \bar{z} = z_0, \ldots, z_{n-1})$ such that for $V_j = \{i < p_j \mid \bar{s}_j \mid : A_j(\bar{x}, \bar{z}, i)\}$

$$\{|z_i| \leq q_i \mid \bar{t}_i \mid \wedge C_i(z_i, \bar{x}) : i < n\}, \Pi \to \Lambda, \{B_j(V_j, \bar{x}) : j < m\}$$

is provable in $AID$.
**Case1** $\Sigma^b_0 - RD$:

$$\exists |y| \leq 2^{\ell \|\bar{t}\|} \forall |i| \leq \ell \|\bar{t}\| [\{0 \neq |i| = \ell \|\bar{t}\| \to (i \in y \leftrightarrow B(\bar{t}, i))\}$$
$$\wedge \{0 \neq |i| < \ell \|\bar{t}\| \to (i \in y \leftrightarrow I(\bar{D}(\bar{t}, i), i0 \in y, i1 \in y))\}]$$

for a sequence $t = \bar{t}(\bar{x})$ of terms. Then put $A(\bar{x}, i) \leftrightarrow_{df} A^{\ell, B, \bar{D}, I}(\bar{t}(\bar{x}), i)$ for the inductively defined $A^{\ell, B, \bar{D}, I}$.
**Case2** Contraction: As in the proof of Lemma 8.1.
**Case3** $\Sigma^b_0 - AC$: By IH we have $i < p \mid t \mid, \Gamma \to \Delta, B(i, \{j < q \mid t \mid : A(\bar{x}, i, j)\})$, where we assume $\Gamma, \Delta \in \Sigma^b_0$ for simplicity. Then $\Gamma \to \Delta, \forall i < p \mid t \mid B(i, Z_i)$ for

$$Z_i = \{k < p \mid t \mid \cdot q \mid t \mid : \exists i < p \mid t \mid \exists j < q \mid t \mid (k = i \cdot q \mid t \mid + j \wedge A(\bar{x}i, j))\}$$

and $Z_i = Z[i \cdot q \mid t \mid, (i+1) \cdot q \mid t \mid)$.
**Case4** Cut: Infer

$$\exists |z| \leq q \mid t \mid C(z, \bar{x}), \Gamma, \Pi \to \Delta, \Lambda, \exists |y| \leq p \mid s \mid B(y, \bar{x})$$

from

$$\exists |z| \leq q \mid t \mid C(z, \bar{x}), \Gamma \to \Delta, \exists |u| \leq r \mid t' \mid D(u, \bar{x})$$

and

$$\exists |u| \leq r \mid t' \mid D(u, \bar{x}), \Pi \to \Lambda, \exists |y| \leq p \mid s \mid B(y, \bar{x})$$

By IH we have

$$|z| \leq q \mid t \mid, C(z, \bar{x}), \Gamma \to \Delta, D(\{i < r \mid t' \mid : A_0(\bar{x}, z, i)\}, \bar{x})$$

and

$$|u| \leq r \mid t' \mid, D(u, \bar{x}), \Pi \to \Lambda, B(\{i < p \mid s \mid : A_1(\bar{x}, u, i)\}, \bar{x})$$

In a proof of the latter sequent, substitute $\{i < r \mid t' \mid : A_0(\bar{x}, z, i)\}$ for the variable $u$ we get

$$D(\{i < r \mid t' \mid : A_0(\bar{x}, z, i)\}, \bar{x}), \Pi \to \Lambda, B(\{i < p \mid s \mid : A_1(\bar{x}, u, i)\}, \bar{x})$$

for $A(\bar{x}, z, i) \Leftrightarrow_{df} A_1(\bar{x}, \{i < r \mid t' \mid : A_0(\bar{x}, z, i)\}, i)$. By a cut with the cut formula $D(\{i < r \mid t' \mid : A_0(\bar{x}, z, i)\}, \bar{x})$ we get

$$|z| \leq q \mid t \mid, C(z, \bar{x}), \Gamma, \Pi \to \Delta, \Lambda, B(\{i < p \mid s \mid : A(\bar{x}, z, i)\}, \bar{x})$$

□

**Theorem 8.1**    1. $\Sigma_0^b - RD + \Sigma_1^b - AC$ is $\Sigma_1^b$-conservative over $AID + \Sigma_0^b - CA$.

2. $\Sigma_0^b - RD + \Sigma_1^b - AC$ is $\Sigma_0^b$-conservative over $AID$.

3. Every $s\Delta_1^b$-formula in $\Sigma_0^b - RD + \Sigma_1^b - AC$ is $\Sigma_0^b$-definable in $AID$: for strict $\Sigma_1^b$-formulae $A, B \in s\Sigma_1^b$, if $\Sigma_0^b - RD + \Sigma_1^b - AC \vdash A \leftrightarrow \neg B$, then $AID \vdash A \leftrightarrow A'$ for a $\Sigma_0^b$ $A'$ in $L_{AID}$.

**Proof**.
8.1.1. Let $C(\bar{x})$ be a $\Sigma_1^b$-formula provable in $\Sigma_0^b - RD + \Sigma_1^b - AC$. Let $B$ be a $\Sigma_0^b$-formula so that for a polynomial $p \, |\bar{x}|$

$$\Sigma_1^b - AC \vdash C(\bar{x}) \leftrightarrow \exists \, |y| \leq p \, |\bar{x}| \, B(y, \bar{x})$$

By Lemma 8.2 and $\Sigma_0^b - CA$, $\exists \, |y| \leq p \, |\bar{x}| \, (y = \{i < p \, |\bar{x}| : A(\bar{x}, i)\})$, we have $AID + \Sigma_0^b - CA \vdash \exists \, |y| \leq p \, |\bar{x}| \, B(y, \bar{x})$. Since $\exists \, |y| \leq p \, |\bar{x}| \, B(y, \bar{x}) \to C(\bar{x})$ is provable without $\Sigma_1^b - AC$, we conclude $AID + \Sigma_0^b - CA \vdash C(\bar{x})$.

8.1.3. Suppose $\Sigma_0^b - RD + \Sigma_1^b - AC \vdash A(\bar{x}) \leftrightarrow \neg B(\bar{x})$ for $A, B \in s\Sigma_1^b$. By Lemma 8.2 there exists a $\Sigma_0^b$-formula $A_0(\bar{x}, i)$ such that

$$AID \vdash (0 \in y \leftrightarrow A_0(\bar{x}, 0)) \to (y = 1 \wedge A(\bar{x})) \vee (y = 0 \wedge B(\bar{x}))$$

By Theorem 8.1.2 we have $B(\bar{x}) \to \neg A(\bar{x})$. Therefore $AID \vdash A_0(\bar{x}, 0) \leftrightarrow A(\bar{x})$.    □

**Corollary 8.1**    1. If a function of polynomial growth rate is $\Sigma_1^b$-definable in $\Sigma_0^b - RD + \Sigma_1^b - AC$ in the sense of [4], then the function is in $\mathcal{F}ALOGTIME$.

2. If a predicate is $\Delta_1^b$-definable in $\Sigma_0^b - RD + \Sigma_1^b - AC$, then the predicate is in $ALOGTIME$.

Corollary 8.1.2, Proposition 1.1 and Theorem 3.1 yields

**Theorem 8.2** *For a predicate $A$,*

$$A \in ALOGTIME \;\;\Leftrightarrow\;\; A \text{ is } \Sigma_0^b\text{-definable in } AID$$
$$\Leftrightarrow\;\; A \text{ is } \Delta_1^b\text{-definable in } \Sigma_0^b - RD + \Sigma_1^b - AC$$

**Theorem 8.3** *For a $\Sigma_0^b$-formula $B(y, \bar{x})$ and a polynomial $p \, |\bar{x}|$, if $\Sigma_0^b - RD + \Sigma_1^b - AC \vdash \forall \bar{x} \exists \, |y| \leq p \, |\bar{x}| \, B(y, \bar{x})$, then there exist $\Sigma_0^b$-formulae $A(\bar{x}, i), P(\bar{x}, i)$ in $L_{AID}$ and a polynomial $q \, |\bar{x}|$ such that*

$$AID \vdash \{i < q \, |\bar{x}| : P(\bar{x}, i)\} \text{ is a Frege proof of } \langle B^*(\{i < p \, |\bar{x}| : A(\bar{x}, i)\}, \bar{x})\rangle,$$

*and hence $AID + \Sigma_0^b - CA \vdash \mathcal{F} \vdash^{q|\bar{x}|} \langle B^*(\{i < p \, |\bar{x}| : A(\bar{x}, i)\}, \bar{x})\rangle$.*

# 9 Clote's $ALV$ and $AID$

In [9] P. Clote defines a function algebra $N_0$ and show that $N_0 = \mathcal{F}ALOGTIME$, the class of $ALOGTIME$-computable functions. Then he introduces an equational system $ALV$ based on $N_0$ in [10].

In this section we show that $AID + \Sigma_0^b - CA$ is equivalent to $ALV$ in the sense that there exist $\Sigma_1^b$-faithful interpretations between $AID + \Sigma_0^b - CA$ and a quantified theory $QALV$.

## 9.1 $AID + \Delta_1^b - CA$ contains $ALV$.

In this subsection we show that (the graph of) each function $f \in N_0$ is $\Delta_1^b$-definable in $AID + \Delta_1^b - CA$. Hence via Clote's result $N_0 = \mathcal{F}ALOGTIME$ in [9] we get the

**Theorem 9.1** *For each $f \in \mathcal{F}ALOGTIME = N_0$ there exists a $\Delta_1^b$-formula $G_f(\bar{x}, y)$ in $AID + \Delta_1^b - CA$ and a polynomial $p_f$ so that $G_f(\bar{x}, y)$ defines the graph of $f$ in the standard model, $AID + \Delta_1^b - CA \vdash \forall \bar{x} \exists ! y G_f(\bar{x}, y)$ and $AID + \Delta_1^b - CA \vdash G_f(\bar{x}, y) \to |y| \leq p_f |\bar{x}|$.*

Further the definitions are intensionaly correct. Namely

**Theorem 9.2** *Using the definition $G_f(i, \bar{x})$ of the graph of $f \in N_0$ the defining equations of $f$ are derivable in $AID + \Delta_1^b - CA$.*

Now we prove Theorems 1 and 2 by induction on the construction of $f \in N_0$.
**Initial functions**. These are zero $o(x) = 0$, successor functions $xi = s_i x = 2 \cdot x + i$ for $i < 2$, projections $i_k^n(x_1, \ldots, x_n) = x_k$, $Bit(i, x)$, $\#$ and the function $tree$. Except the last one $tree$, which is $NC^1$-complete for $AC^0$ reduction, cf. [9], the assertions are clear.

The function $tree$ takes values $0, 1$ and so we regard it as a predicate. Then the predicate $tree$ is defined from the auxiliary functions $and(x), or(x)$ as follows: First set $and(0) = or(0) = 0 \ \& \ and(x) = or(x) = 1$ for $1 \leq x \leq 3$. For $x > 0$

$$\begin{array}{ll} and(x00) = and(x)0 & or(x00) = or(x)0 \\ and(x10) = and(x)0 & or(x10) = or(x)1 \\ and(x01) = and(x)0 & or(x01) = or(x)1 \\ and(x11) = and(x)1 & or(x11) = or(x)1 \end{array}$$

Then
$$tree(x) \Leftrightarrow \begin{cases} parity(x) =_{df} Bit(0, x) = 1 & \text{if } x < 16 \\ tree(or(and(x))) & \text{otherwise} \end{cases}$$

Therefore, for $x > 1$, $tree(x)$ is the predicate obtained by evaluating a perfect $and/or$ tree on the $4^{\lfloor \log_4(|x| \dot{-} 1) \rfloor}$ many least significant bits of $x$. We define a predicate $Tree(x, p)$ inductively so that $tree(x) \leftrightarrow Tree(x, 1)$. Put

$$y = x[0, 4^{\lfloor \log_4(|x| \dot{-} 1) \rfloor})$$

for $x > 1$ and put $y = 0$ for $x \leq 1$. If $x > 1$, then $|y| = 4^{\lfloor \log_4(|x| \dot{-} 1) \rfloor} = \max\{u < |x| : \exists z(u = 4^z)\}$ and hence

$$\|y\| \dot{-} 1 = \begin{cases} \lfloor \frac{\|x\| \dot{-} 2}{2} \rfloor & \text{if } \exists k(|x| = 2^k)[\leftrightarrow \forall i < \|x\| \dot{-} 1 (Bit(i, |x|) = 0)] \\ \lfloor \frac{\|x\| \dot{-} 1}{2} \rfloor & \text{otherwise} \end{cases} \quad (10)$$

Also $\|y\| \dot{-} 1 \leq \|x\|$.

**(T.0)** $Tree(x, p) \to 0 \neq |p| \leq \|x\| + 1$.

**(T.1)** The case $\|y\| \leq |1p| \leq \|x\| + 1$:

$$Tree(x, 1p) \leftrightarrow Bit(p[0, \|y\| \dot{-} 1], x) = 1$$

**(T.2)** The case $0 \neq |p| < \|y\|$:

$$Tree(x, p) \leftrightarrow [\{|p| \text{ is odd } \& (Tree(x, p0) \lor Tree(x, p1))\}$$
$$\lor \{|p| \text{ is even } \& (Tree(x, p0) \land Tree(x, p1))\}]$$

Note that by (10) this defines the predicate $Tree = A^{\ell, B, \bar{D}, I}$ in $L_{AID}$ for some $\Sigma_0^b$-formulae $B, \bar{D}$ in $L_{BA}$ and a boolean $I$ with $\ell\|x\| = \|x\| + 1$.

Now we show $Tree(x, 1)$ enjoys the defining axioms of $tree$ in $AID$. First assume $x < 16$. Then $|x| \leq 4$ and hence $|y| = 1 \,\&\, \|y\| \dot{-} 1 = 0$. By **(T.1)** with $p = 0$ we have $Tree(x, 1) \leftrightarrow Bit(0, x) = parity(x) = 1$.

Next consider the case $x \geq 16$. We have to show

$$Tree(x, 1) \leftrightarrow Tree(or(and(x)), 1)$$

We understand this formula is an abbreviation for the $\Delta_1^b$-formula

$$\forall |y| \leq |x| [y = or(and(x)) \to (Tree(x, 1) \leftrightarrow Tree(y, 1))] \leftrightarrow$$
$$\exists |y| \leq |x| [y = or(and(x)) \,\&\, (Tree(x, 1) \leftrightarrow Tree(y, 1))]$$

for a $\Sigma_0^b$-formula $y = or(and(x))$.

Observe that $|or(x)| = |and(x)| = \lfloor \frac{|x|}{2} \rfloor + parity(|x|)$ for $x > 3$. We show the follwing Claim.

**Claim 9.1** *Put*

$$z = 2\lfloor \log_4(|or(and(x))| \dot{-} 1) \rfloor = 2\lfloor \log_4(|x| \dot{-} 1) \rfloor - 2$$

*For $x \geq 16$ and any odd $|1p| \leq z + 1$*

$$Tree(x, 1p) \leftrightarrow Tree(or(and(x)), 1p)$$

**Proof** of Claim 9.1 by induction on $z+1-|1p|$. Here we use $\Delta_1^b - LIND$. This follows from $\Delta_1^b - CA$. If $|1p| < z+1$, then the Claim follows from IH and **(T.2)**. Suppose $z+1 = |1p|$. Then by **(T.1)** we have $Tree(or(and(x)), 1p) \leftrightarrow Bit(p, or(and(x))) = 1$. On the other we have by **(T.2)** and **(T.1)**

$$Tree(x, 1p) \leftrightarrow$$
$$(Tree(x, 1p00) \wedge Tree(x, 1p01)) \vee (Tree(x, 1p10) \wedge Tree(x, 1p11)) \leftrightarrow$$
$$(Bit(p00, x) = 1 \wedge Bit(p01, x) = 1) \vee (Bit(p10, x) = 1 \wedge Bit(p11, x) = 1)$$

Put $n = 4k + j = |x| \geq 5$ for $j < 4$. We show

$$[j \neq 0 \Rightarrow 4p + 3 < 4k] \, \& \, [j = 0 \Rightarrow 4p + 3 < 4k - 4] \tag{11}$$

By $|1p| = 2\lfloor \log_4(n \dot{-} 1) \rfloor - 1$ we have $|p| \leq 2\lfloor \log_4(n \dot{-} 1) \rfloor - 2$ and hence $p + 1 \leq 2^{2\lfloor \log_4(n \dot{-} 1) \rfloor - 2}$. Therefore $4p + 4 \leq 2^{2\lfloor \log_4(n \dot{-} 1) \rfloor} = 4^{\lfloor \log_4(n \dot{-} 1) \rfloor}$.

**Case1**. $j \neq 0$: Then $4^{\lfloor \log_4(n \dot{-} 1) \rfloor} = 4^{\lfloor \log_4 4k \rfloor} \leq 4k$. Hence $4p + 3 < 4k$.

**Case2**. $j = 0$: Then $k \geq 2$ by $n = 4k \geq 5$. Put $n - 1 = 4k - 1 = \sum_{i=0}^{m} 4^i \cdot y_i$ with $0 \leq y_i < 4$, $y_m \neq 0$ & $y_0 = 3$. Then $4^{\lfloor \log_4(n \dot{-} 1) \rfloor} = 4^m$. It suffices to show $4k - 4 - 4^m \geq 0$ to have $4p + 3 < 4k - 4$. We have $4k - 4 - 4^m = 4^m \cdot (y_m - 1) + \sum_{i < m} 4^i \cdot y_i - 3$. Suppose $m = 0$. Then we would have $4 \leq n - 1 = 4k - 1 \leq 3$. A contradiction. Hence $m > 0$. The assertion follows from $m > 0$ and $y_0 = 3$.

Thus we have shown (11) and from this and $4p + 3 = p11$ we see

$$Bit(p, or(and(x))) = 1 \leftrightarrow Bit(p0, and(x)) = 1 \vee Bit(p1, and(x)) = 1 \leftrightarrow$$
$$(Bit(p00, x) = 1 \wedge Bit(p01, x) = 1) \vee (Bit(p10, x) = 1 \wedge Bit(p11, x) = 1)$$

Thus we have $Tree(x, 1p) \leftrightarrow Tree(or(and(x)), 1p)$ as desired.

*End of Proof of Claim*

**Constructors**. These are compositions and *CRN* (Concatenation Recursion on Notation): define $f$ from $g$ and $h_i$ $(i < 2)$ by

$$f(0, \bar{x}) = g(\bar{x}) \,;\, n > 0 \to f(n, \bar{x}) = s_i(f(\lfloor \frac{n}{2} \rfloor, \bar{x}))$$

where $i = sg(h_j(\lfloor \frac{n}{2} \rfloor, \bar{x}))$ with $j = parity(n) = Bit(0, n)$ and $sg(x) < 2$ & $(sg(x) = 0 \leftrightarrow x = 0)$. Composition is harmless and for *CRN* define

$$p_f(|n|, |\bar{x}|) = p_g |\bar{x}| + |n|$$

and

$$G_f(n, \bar{x}, y) \leftrightarrow \exists |z| \leq p_g |\bar{x}| [g(\bar{x}) = z \,\&\, \exists |u| \leq |n| \{y = u * z' \,\&$$
$$\forall i < |n| \, (Bit(i, u) = 0 \leftrightarrow \bigvee_{j<2} (h_j(\lfloor \frac{n}{2} \rfloor, \bar{x}) = 0 \,\&\, Bit(i, n) = j))\}]$$

where $z' = 2^{|z|} + z$, $g(\bar{x}) = z$ denotes $G_g(\bar{x}, z)$ and $h_j(n, \bar{x}) = i$ denotes $G_h(n, \bar{x}, i)$. Note that by IH $G_h(n, \bar{x}, i)$ is $\Delta_1^b$ in $AID + \Delta_1^b - CA$ and hence the RHS is $\Sigma_1^b$. Further a $\Pi_1^b$-form of $G_f(n, \bar{x}, y)$ is defined similarly.

To show $\exists y G_f(n, \bar{x}, y)$ it suffices to show

$$\exists |u| \leq |n| \, \forall i < |n| \, (Bit(i, u) = 0 \leftrightarrow \bigvee_{j<2} (h_j(\lfloor \frac{n}{2} \rfloor, \bar{x}) = 0 \,\&\, Bit(i, n) = j))$$

and this follows from $\Delta_1^b - CA$.

Thus we have shown Theorems 1 and 2.

## 9.2 Inductive definition in $ALV$

In this subsection we show that $ALV$ can simulate inductive definitions in $AID$. Specifically it is shown that for each inductively defined predicate $A$ in $AID$, there exists a $\{0,1\}$-valued function symbol $f_A$ in $ALV$ such that $f_A(\bar{x}, p) = 1$ satisfies the definig axioms **(A.0)-(A.2)** of the predicate $A$ demonstrably in $ALV$.

Let $QALV$ denote a quantified version of $ALV$. It is a first order theory whose non-logical symbols are those of $ALV$ and whose axioms are the universal closures of defining equations of these function symbols and induction on notation, together with two more axioms: $0 \neq 1$ and $\forall x (\lfloor \frac{x}{2} \rfloor = 0 \rightarrow (x = 0 \vee x = 1))$. The last two axioms are added by Cook [15]. Note that in [15] the same name $QALV$ designates a different first order theory, i.e., a quantified version of $ALV'$ in [11]. As in [15] we see easily that

**Proposition 9.1**  1. $QALV$ is a conservative extension of $ALV$.

2. For each $\Sigma_0^b$-formula $B$ in $QALV$ there exists a function symbol $f_B$ such that $QALV \vdash B(\bar{x}) \leftrightarrow f_B(\bar{x}) = 0$, and hence using $CRN$ we have $QALV \vdash \Sigma_0^b - CA$: for each $\Sigma_0^b$-formula $B$ and each polynomial $p$, $QALV \vdash \forall \bar{x} \exists y (y = \{i < p \,|\, \bar{x} |: B(\bar{x}, i)\})$.

3. $QALV$ proves $\Sigma_0^b - LIND$.

**Theorem 9.3** *For each inductively defined predicate $A$ in $AID$, there exists a $\{0,1\}$-valued function symbol $f_A$ in $ALV$ such that $f_A(\bar{x}, p) = 1$ satisfies the definig axioms* **(A.0)-(A.2)** *of the predicate $A$ demonstrably in $QALV$.*

This together with Proposition 9.1 yields the

**Corollary 9.1** $QALV$ *and hence* $ALV$ *proves reflection schema* $RFN(\mathcal{F})$ *for a Frege system.*

Now we prove Theorem 9.3. Let $A = A^{\ell, B, \bar{D}, I}$ be a given inductively defined predicate in $AID$. Recall that for a formula $F$ and $i < 2$ we have putted

$$F^i = \begin{cases} F & i = 1 \\ \neg F & i = 0 \end{cases} \tag{2}$$

For $i < 2$ put $k(i) = Bit(i, k)$.

First convert the boolean formulae $I$ and $\neg I$ into DNF to yield: for $\xi < 2$,

$$I^\xi(\bar{D}(x, p), A(\bar{x}, p0), A(\bar{x}, p1)) \leftrightarrow$$
$$\bigvee \{I_k^\xi(\bar{x}, p) \wedge \bigwedge \{A(\bar{x}, pi)^{k(i)} : i < 2\} : k < 2^2\} \leftrightarrow$$
$$[\bigvee \{I_k^\xi(\bar{x}, p) \wedge \bigwedge \{A(\bar{x}, pi)^{k(i)} : i < 2\} : k = 0, 1\} \wedge (\top \vee \top)] \vee$$
$$[\bigvee \{I_k^\xi(\bar{x}, p) \wedge \bigwedge \{A(\bar{x}, pi)^{k(i)} : i < 2\} : k = 2, 3\} \wedge (\top \vee \top)] \tag{12}$$

Here is an *and/or* tree of depth 4. We define function symbols $g_A^\xi$ so that for $\xi, \eta < 2$

$$tree(g_A^\xi(\bar{x}, p)) = \eta \Leftrightarrow A^\lambda(\bar{x}, p)$$

where $\lambda = (\xi \leftrightarrow \eta)$.

The resulting perfect *and/or* tree is of depth $4(\ell\|\bar{x}\| - |p|)$.

In view of **(A.0)**, put

$$\neg(0 \neq |p| \leq \ell\|\bar{x}\|) \to g_A^1(\bar{x}, p) = 0 \ \& \ g_A^0(\bar{x}, p) = 1$$

In what follows suppose $0 \neq |p| \leq \ell\|\bar{x}\|$. We define function symbols $h_A^\xi$ so that $g_A^\xi(\bar{x}, p) = 2^{2^{4(\ell\|\bar{x}\| - |p|)}} + h_A^\xi(\bar{x}, p)$ do the job. Put $y^\xi = h_A^\xi(\bar{x}, p)$.

First consider the case $|p| = \ell\|\bar{x}\|$. Then by **(A.1)** we have $A(\bar{x}, p) \leftrightarrow B(\bar{x}, p)$. Put, *cf.* Proposition 9.1.2,

$$y^\xi = \begin{cases} 1 & \text{if } B^\xi(\bar{x}, p) \\ 0 & \text{otherwise} \end{cases}$$

Next assume $|p| < \ell\|\bar{x}\|$. In the following we give a $\Sigma_0^b$ condition in $L_{BA}$ which is equivalent to $Bit(q, y^\xi) = 1$ for $q < |y^\xi| \leq 2^{4(\ell\|\bar{x}\| - |p|)}$ and hence $|q| \leq 4(\ell\|\bar{x}\| - |p|)$. Then by $CRN$ we can pick a desired function symbol $h_A^\xi$.

**Case1**. $\exists i \leq \lfloor \frac{|q|}{4} \rfloor (Bit(4i + 2, q) = 1)$: Put $Bit(q, y^\xi) = 1$. This corresponds to $(\top \vee \top)$ in (12). Namely once the branch $q$ has entered in definitely true subtree corresponding to $(\top \vee \top)$, the leaf of $q$ receives 1.

**Case2**. Suppose $\neg[\exists i \leq \lfloor \frac{|q|}{4} \rfloor (Bit(4i + 2, q) = 1)]$. For $m \leq \ell\|\bar{x}\| - |p|$ put

$$p_m = \{i < m : Bit(4i, q[i_1 - 4m, i_1)) = 1\}$$

with $i_1 =_{df} 4(\ell\|\bar{x}\| - |p|)$. Note that $q[i_1 - 4m, i_1)$ denotes the prefix of length $4m$ of the string $0^{[i_1 - |q|]} * q$ over $\{0, 1\}$ of length $i_1$.

Put $i_0 =_{df} i_1 - 4m$. We say that $q$ is *positive at* $m > 0$ if

$$(Bit(i_0, q) = 0 \,\&\, Bit(i_0 + 1, q) = 1) \lor (Bit(i_0, q) = 1 \,\&\, Bit(i_0 + 3, q) = 1)$$

This means that the subtree determined by the prefix $q[i_0, i_1)$ evaluates the value of $A(\bar{x}, p \cdot 2^m + p_m)$, cf (12). Also we say that $q$ is *positive at* $m = 0$ if $\xi = 1$. Otherwise $q$ is said to be *negative at* $m \geq 0$.

Put $j = i_0 - 4$ and let $k = Bit(j+3, q)Bit(j+1, q)$ denote $k < 4 \,\&\, Bit(1, k) = Bit(j+3, q) \,\&\, Bit(0, k) = Bit(j+1, q)$ if $j = i_0 - 4 = i_1 - 4m - 4 \geq 0$, i.e., if $m < \ell\|\bar{x}\| - |p|$. Otherwise, viz. in the case $m = \ell\|\bar{x}\| - |p|$, let $k = Bit(j+3, q)Bit(j+1, q)$ denote a $k < 4$ and put $I_k^\xi(\bar{x}, p \cdot 2^m + p_m) \Leftrightarrow_{df} B^\xi(\bar{x}, p \cdot 2^m + p_m)$ for any $k < 4$.

Now the following is the necessary and sufficient condition to be $Bit(q, y^\xi) = 1$ in this **Case2**:

$$\forall m \leq \ell\|\bar{x}\| - |p| \bigwedge_{k<4} [k = Bit(j+3, q)Bit(j+1, q) \rightarrow$$
$$\{q \text{ is positive at } m \Rightarrow I_k^1(\bar{x}, p \cdot 2^m + p_m)\} \,\&\,$$
$$\{q \text{ is negative at } m \Rightarrow I_k^0(\bar{x}, p \cdot 2^m + p_m)\}$$

These give a $\Sigma_0^b$-definition of $Bit(q, y^\xi) = 1$. Then $\{(\bar{x}, p) : tree(g_A^\xi(\bar{x}, p)) = 1\}$ satisfies the defining axioms **(A.0)-(A.2)** for $A = A^{\ell, B, \bar{D}, I}$. Lemma 9.1.2, 9.1.3 yields Theorem 9.3.

**Lemma 9.1** *1. Suppose $z = 2^{2^{4m}} + y$ with $m > 0$ and $y < 2^{2^{4m}}$. For $k < 2^4$ put $z_k = 2^{2^{4(m-1)}} + y[2^{4(m-1)}k, 2^{4(m-1)}(k+1))$. Let $y_0 < 2^{2^4}$ denote the number such that $Bit(k, y_0) = tree(z_k)$ for $k < 2^4$. Then*

$$tree(z) = 1 \leftrightarrow or\bigl(and\bigl(or\bigl(and(2^{2^4} + y_0)\bigr)\bigr)\bigr) = 1$$

*2. Suppose $0 \neq |p| < \ell\|\bar{x}\|$. Then for $i_1 = 4(\ell\|\bar{x}\| - |p|)$*

$$tree(2^{2^{i_1}} + h_A^\xi(\bar{x}, p)) = 1 \leftrightarrow$$
$$\bigvee\{I_k^\xi(\bar{x}, p) \land \bigwedge\{tree(2^{2^{i_1-4}} + h_A^{k(i)}(\bar{x}, pi)) = 1 : i < 2\} : k < 2^2\}$$

*3. Suppose $0 \neq |p| \leq \ell\|\bar{x}\|$. Then*

$$tree(2^{2^{4(\ell\|\bar{x}\|-|p|)}} + h_A^1(\bar{x}, p)) = 1 \Leftrightarrow tree(2^{2^{4(\ell\|\bar{x}\|-|p|)}} + h_A^0(\bar{x}, p)) = 0$$

**Proof**.
9.1.1. This is proved by induction on $m > 0$.
9.1.2. Use Lemma 9.1.1.
9.1.3. This is proved by induction on $\ell\|\bar{x}\| - |p|$ using Lemma 9.1.2. □

## 9.3 $\Sigma_1^b$-faithful interpretations

In this subsection we conclude that translations derived from the prceeding subsections are $\Sigma_1^b$-faithful.

For a formula $B$ in $QALV$ let $I_{ID}(B)$ denote the formula in $AID$ which is obtained from $B$ by replacing each $f \in N_0$ by its $\Delta_1^b$-graph defined in the proof of Theorems 9.1 and 9.2. Observe that for $\Sigma_1^b$-formula $B$ in $QALV$ $I_{ID}(B)$ is a $\Sigma_1^b$-formula in $AID$.

For a formula $B$ in $AID$ let $I_V(B)$ denote the formula in $QALV$ which is obtained from $B$ by replacing each $A^{\ell,B,\bar{D},I}(\bar{x},p)$ by a $tree(g_A^\xi(\bar{x},p)) = 1$ defined in the proof of Theorems 9.3. Observe that for $\Sigma_1^b$-formula $B$ in $AID$ $I_V(B)$ is a $\Sigma_1^b$-formula in $QALV$.

**Theorem 9.4** 1. For each $\Sigma_1^b$-formula $B$ in $QALV$,

$$QALV \vdash B \Leftrightarrow AID + \Sigma_0^b - CA \vdash I_{ID}(B)$$

2. For each $\Sigma_1^b$-formula $B$ in $AID$,

$$AID + \Sigma_0^b - CA \vdash B \Leftrightarrow QALV \vdash I_V(B)$$

**Proof**. Let $T$ denote temporalily a union of theories $AID + \Sigma_0^b - CA$ and $QALV$: its language $L_T$ is the union $L_{AID} \cup L_{QALV}$, and its axioms are ones in $AID + \Sigma_0^b(L_{AID}) - CA$ and $QALV$ plus $\Sigma_0^b(L_T) - LIND$. By Theorems 9.2 and 9.3 we have

1. For $\Sigma_1^b$-formula $B$ in $QALV$, $T + \Delta_1^b(L_{AID}) - CA \vdash B \leftrightarrow I_{ID}(B)$.

2. For $\Sigma_1^b$-formula $B$ in $AID$, $T \vdash B \leftrightarrow I_V(B)$.

9.4.1. For a $\Sigma_1^b$-formula $B$ in $QALV$ first suppose $QALV \vdash B$. Then $T + \Delta_1^b(L_{AID}) - CA \vdash I_{ID}(B)$. By replacing any formula $C$ in this proof by $I_{ID}(C)$ (leave formulae in $AID$ unchanged) we get $AID + \Delta_1^b(L_{AID}) - CA \vdash I_{ID}(B)$. Thus by Theorem 8.1.1 we conclude $AID + \Sigma_0^b - CA \vdash I_{ID}(B)$.

Conversely assume $AID + \Sigma_0^b - CA \vdash I_{ID}(B)$. Then $T + \Delta_1^b(L_{AID}) - CA \vdash B$. By replacing any formula $C$ in this proof by $I_V(C)$ (leave formulae in $QALV$ unchanged) we get $QALV + \Delta_1^b(L_{QALV}) - CA \vdash B$. As in the proof of Theorem 8.1.1 or cf. Theorem 6 in p.79, [15], we see that $QALV + \Sigma_1^b(L_{QALV}) - AC$ and hence $QALV + \Delta_1^b(L_{QALV}) - CA$ is $\Sigma_1^b$-conservative over $QALV$. Thus we conclude $QALV \vdash B$.

9.4.2. This is seen as in Theorem 9.4.1 using Proposition 9.1. □

# References

[1] M. Ajtai, The complexity of the pigeonole principle, in Proceedings of the $29^{th}$ Annual IEEE Symposium on Foundations of Computer Science, 1988, 346-355.

[2] J. L. Balcázar, J. Díaz and J. Gabarró, Structural Complexity II, EATCS Monographs on Theoretical Computer Sci. 22, (Springer, Berlin, 1988)

[3] M. Bonet, S. Buss and T. Pitassi, Are there hard examples for Frege systems?, in:P. Clote and J. B. Remmel, eds., Feasible Mathematics II (Birkhäuser, Boston Basel Berlin, 1995), 30-56.

[4] S. R. Buss, Bounded Arithmetic (Bibliopolis, Napoli, 1986).

[5] S. R. Buss, Polynomial size proofs of the propositional pigeonhole principles, J Symb Logic 52(1987), 916-927.

[6] S. R. Buss, The Boolean formula value problem is in $ALOGTIME$, Proceedings of the $19^{th}$ Annual ACM Symposium on Theory of Computing, May 1987, 123-131.

[7] S. R. Buss, Propositional consistency proofs, Ann Pure Appl Logic 52(1991), 3-29.

[8] S. R. Buss, Algorithms for Boolean formula evaluation and tree contraction, in Arithmetic, Proof Theory and Computational Cimplexity, P. Clote and J. Krajíček, eds., Oxford UP, 193, 96-115.

[9] P. Clote, Sequential, machine-independent characterizations of the parallel complexity classes $ALOGTIME$, $AC^k$, $NC^k$ and $NC$, in: S.R. Buss and P. Scott, eds., Feasible Mathematics (Birkhäuser, Boston Basel Berlin, 1990) 49-70.

[10] P. Clote, $ALOGTIME$ and a conjecture of S. A. Cook, Annals of Mathematics and Artificial Intelligence 6(1992), 57-106. extended abstract in Proceedings of IEEE Logic in Computer Science, Philadelphia, June 1990.

[11] P. Clote, On polynomial size Frege proofs of certain combinatorial pinciples, in: P. Clote and J. Krajíček, eds., Arithmetic, Proof Theory and Computational Complexity (Oxford UP, Oxford, 1993) 162-184.

[12] P. Clote and G. Takeuti, Bounded Arithmetics for $NC$, $ALOGTIME$, $L$ and $NL$, Ann Pure Appl Logic 56(1992), 73-117.

[13] P. Clote and G. Takeuti, First order bounded arithmetic and small boolean circuit complexity classes, in:P. Clote and J. B. Remmel, eds., Feasible Mathematics II (Birkhäuser, Boston Basel Berlin, 1995) 154-218

[14] S. A. Cook, Feasibly constructive proofs and the propositional calculus, in Proceedings of the $7^{th}$ Annual ACM Symposium on Theory of Computing, 1975, 83-97.

[15] S. A. Cook, Relating the provable collapse of $P$ to $NC^1$ and the power of logical theories, in DIMACS series in Discrete Mathematics and Theoreical Computer Science Volume 39 (1998), 73-91

[16] J. Krajíček, On Frege and extended Frege proof systems, in:P. Clote and J. B. Remmel, eds., Feasible Mathematics II (Birkhäuser, Boston Basel Berlin, 1995), 284-319.

[17] J. Krajíček, Bounded arithmetic, propositional logic, and complexity theory, (Cambridge UP, Cambridge, 1995)

[18] F. Pitt, A quantifier-free theory based on a string algebra for $NC^1$, in DIMACS series in Discrete Mathematics and Theoreical Computer Science Volume 39 (1998),229-252.